\input amstex
\input epsf
\documentstyle{amsppt}
\NoBlackBoxes
\magnification=1200
\parindent 20 pt
\vsize=7.50in
\define\ve{\varepsilon}
\define\1{^{-1}}
\define\s{\sigma}
\define\be{\beta}
\define\CP{\Bbb C \Bbb P}
\define\F{\Bbb F _1}
\define\N{\Bbb N}

\define \tg{\tilde{\gamma}}

\define \G{\Gamma}

\define \M{\Cal M}
\define \vp{\varphi}
\define \Dl{\Delta}
\define \dl{\delta}
\define \C{\Bbb C}

\define \p{\partial}
\define \ri{\rightarrow}
\define \usr{\underset\sim\to\rightarrow}
\define \Int{\operatorname{Int}}
\define \Center{\operatorname{Center}}

\define \mk{\medskip}
\define \tm{\times}
\define \bk{\bigskip}
\define \ov{\overline}
\define \edm{\enddemo}
\define \ep{\endproclaim}

\topmatter
\title Braid Monodromy Factorization and Diffeomorphism Types\endtitle
\author Vik. S. Kulikov and M. Teicher\endauthor
\abstract In this paper we prove that if two   cuspidal plane curves $B_1$
and $B_2$ have
equivalent braid monodromy factorizations, then $B_1$ and $B_2$ are
smoothly isotopic in
$\CP^2$.
As a consequence, we obtain that if the discriminant curves (or branch curves
in other terminology) $B_1$ and $B_2$
of generic
projections to $\CP^2$ of surfaces of general type $S_1$ and $S_2,$
imbedded in a projective
space by means of a multiple canonical class have equivalent braid
monodromy factorizations,
then $S_1$ and $S_2$ are diffeomorphic (if we consider them as real
4-folds).\endabstract
\thanks \flushpar Partially supported by INTAS-OPEN-972072.\newline
Partially supported by the Max-Planck Institut, Bonn.\newline
Partially supported by the Emmy Noether Research Institute for Mathematics,
Bar-Ilan
University and the Minerva Foundation, Germany.\endthanks
\endtopmatter

\document

\baselineskip 20pt

\subheading{\S 0.\ Introduction}

Let $S$ be a non-singular algebraic surface in a projective space $\CP ^r$ of
$\deg S=N$. It is well-known that for almost all projections
$pr:\CP ^r \to \CP ^2$ the restrictions $f:S\to \CP ^2$ of
these projections to $S$ satisfy the following four conditions:
\roster
\item"(i)" $f$ is a finite morphism of $\deg f=\deg S$;
\item"(ii)" $f$ is branched along an irreducible curve $B\subset \CP ^2$
with ordinary cusps and nodes, as the only singularities;
\item"(iii)"  $f^{*}(B)=2R+C$, where $R$ is irreducible and non-singular,
and $C$
is reduced;
\item"(iv)"  $f_{\mid R}:R\to B$ coincides with the
normalization of $B$.
\endroster
We will call such $f$ {\it a generic morphism} and its branch curve will be
called {\it the discriminant curve} of $f$.

Two generic morphisms $(S_1,f_1)$, $(S_2,f_2)$ with the same discriminant
curve $B$ are said to be equivalent if there exists an isomorphism
$\varphi : S_1 \to S_2$ such that
$f_1=f_2\circ \varphi $.

The following assertion is known as Chisini's Conjecture.
\proclaim{Chisini's Conjecture}
Let $B$ be the discriminant curve of a generic morphism $f:S\to \CP ^2$
of $\deg f \geq 5$. Then $f$ is uniquely determined by the pair
$(\CP ^2,B)$.
\endproclaim

It is easy to see that the similar conjecture for generic morphisms of
projective curves to $\CP ^1$ is not true. On the other hand in
\cite{Kul1}
it is shown that Chisini's Conjecture holds for the discriminant curves of
almost all generic morphisms of any projective surface. In particular,
if $S$ is any surface of general type with  ample canonical class, then
Chisini's Conjecture
holds for the discriminant curves of the generic morphisms $f:S\to \CP ^2$
given
by a three-dimensional linear subsystem of the $m$-canonical class of
$S$, where
$m\in \N$. The discriminant curve of such a  generic morphism will be
called {\it
$m$-canonical discriminant curve}.

Let $B$ be an algebraic curve in
$\CP ^2$ of degree $p$. Topology of the embedding  $B\subset \CP ^2$ is
determined by the {\it braid monodromy} of $B$
which is described by a factorization of the ``full-twist"
$\Delta _{p}^2$ in the semi-group $B^{+}_{p}$ of the braid group
$B_{p}$ of $p$
string braids (in standard generators, $\Delta _{p}^2=
(X_1\cdot ...\cdot X_{p-1})^{p}$). If $B$ is a cuspidal curve,
then this factorization can be written as follows
$$
\Delta _{p}^2=\prod_{i}Q_i^{-1}X_1^{\rho _i}Q_i, \, \, \, \, \, \, \, \, \, \,
\rho _i\in (1,2,3), \tag 1
$$
where $X_1$ is a positive half-twist in $B_{p}$.

Let
$$
h=g_1\cdot ...\cdot g_r \tag 2
$$
be a factorization in $B^{+}_{p}$. The transformation which changes two
neighbouring factors in (2) as follows:
$$g_i\cdot g_{i+1} \longmapsto (g_ig_{i+1}g_i^{-1})\cdot g_i, $$
or
$$g_i\cdot g_{i+1} \longmapsto g_{i+1}(g_{i+1}^{-1}g_ig_{i+1}) $$
is called a {\it Hurwitz move}.

For $z\in B_{p}$, we denote
$$h_z=z^{-1}g_1z\cdot z^{-1}g_2z\cdot ...\cdot z^{-1}g_rz$$
and say that the factorization expression $h_z$ is obtained from (2)
by simultaneous conjugation by $z$. Two factorizations are called
{\it Hurwitz and conjugation equivalent} if one can be obtained from the
other by a finite sequence of Hurwitz moves followed by a simultaneous
conjugation.
We will say that two factorizations of the form (1) belong to the
same {\it braid factorization type} if they are Hurwitz and conjugation
equivalent. The main problems in this direction are the following:
\proclaim{Problem 1}
Let $B\subset \CP ^2$ be a cuspidal curve. Does the braid factorization type
of the pair $(\CP ^2,B)$ uniquely
determine the diffeomorphic type of this pair, and vice
versa?
\endproclaim
\proclaim{Problem 2}
Let $\Delta^2_{p}={\Cal{E}}_1$ and $\Delta^2_{p}={\Cal{E}}_2$ be two braid
monodromy factorizations. Does there exist a finite algorithm to
recognize whether these two braid monodromy factorizations belong to the
same braid factorization type?
\endproclaim

One of the main results of this article is
\proclaim{Theorem 1}
Let $B_1,\, B_2\subset \CP ^2$ be two cuspidal algebraic curves. Assume that
the pairs $(\CP ^2,B_1)$ and $(\CP ^2,B_2)$ have the same braid factorization
type. Then the pairs $(\CP ^2,B_1)$ and $(\CP ^2,B_2)$ are diffeomorphic.
\endproclaim

It is well-known that there exist four-dimensional smooth manifolds which are
homeomorphic, but not diffeomorphic.
One of the most important problems in four-dimensio\-nal geometry is to
find invariants which distinguish smooth
structures on the
same topological four-dimensional manifold. We believe that in the
algebraic case one can use
the braid factorization type of the discriminant curve of a generic
morphism of a projective surface $S$ to $\CP ^2$ as an invariant of
the smooth structure
(induced by complex structure) on $S$,
considered as a four-dimensional real manifold.
We will prove
\proclaim{Theorem 2}
Let $f_1:S_1\to \CP ^2$ and $f_2:S_2\to \CP ^2$ be two generic morphisms of
non-singular projective surfaces, and let $B_1,\, B_2\subset \CP ^2$ be
their discriminant curves. Assume that Chisini's Conjecture holds for
$(\CP ^2, B_1)$. If the pairs $(\CP ^2,B_1)$ and $(\CP ^2,B_2)$ have the same
braid factorization type, then $S_1$ and $S_2$ are diffeomorphic.
\endproclaim
\proclaim{Corollary}
Let $S_1$ and $S_2$ be two surfaces of general type with ample canonical
class and let $B_1$ and $B_2$ be $m$-canonical discriminant curves,
respectively, of
generic morphisms $f_1:S_1\to \CP ^2$ and $f_2:S_2\to \CP ^2$ given
by three-dimensional linear subsystem of $m$-canonical class on $S_i$, where
$m\in \N$. If the pairs $(\CP ^2,B_1)$ and $(\CP ^2,B_2)$ have the same
braid factorization type, then $S_1$ and $S_2$ are diffeomorphic.
\endproclaim

Here is a brief summary of the rest of this article.
In Sections 1 - 4, we recall definitions and some facts on the braid monodromy
technique developed by B. Moishezon and the second author.
Section 5 is devoted to the description of generators of the centralizer of
the multiple half-twists in the braid group. This description is a key to
the proof of Theorem 1. In Section 6, we recall or prove
some assertions (which may be well-known) on
smooth isotopy of manifolds which we will use in the proof of the main
results. Section 7 is devoted to the proof of Theorem 1, and in Section 8 we
prove Theorem 2.

\bk

\subheading{\S 1.\ Braid monodromy of an affine curve}

Throughout this paper, we  use the following notations:

$S$ is a curve in $\C^2,$ \ $p=\deg S$,

$\pi: \C^2\ri \C$ a generic projection on the first coordinate,

$K(x)=\{y\bigm| (x,y)\in S\}$ \ $(K(x)=$ projection to $y$-axis of
$\pi^{-1}(x)\cap S)$,

$N=\{x\bigm| \# K(x)\lvertneqq p\}$,

$M'=\{x\in S\bigm| \pi\bigm |_X$ is not \'etale at $x\}\ (\pi(M')=N)$.

Assume $\#(\pi\1(x)\cap M^1)=1$ for all $x\in N$.

Let $E$ (resp. $D$) be a closed disk on $x$-axis (resp. $y$-axis) s.t.
$M'\subset E\times D$, $N\subset\Int(E)$, and such that
$\pi _{|(E\times D)\cap B}$ is a proper morphism of degree $p$.

We choose $u\in\p E$ and put
$K=K(u)=\{q_1,\dots,q_p\}.$

In such a situation, we will introduce the concept of ``braid monodromy."

\definition{Definition} \ \underbar{Braid monodromy of $B$ with respect to
$E\times D$, $\pi $, $u$.}

Every loop $\gamma :[0,1]\to E\setminus N$ starting at $u$ has liftings
to a system of $p$ paths in $(E\setminus N)\times D$ starting at
$q_1,\dots, q_p.$
Projecting them to $D$ we get $p$ paths in $D$ defining a motion
$\{q_1(t),\dots,q_p(t)\}$ of $p$ points in $D$ starting and ending at $K.$

This motion defines a braid in $B_p[D,K],$ (see \cite{MoTe1}, Section III).
Thus we get a map $\vp:\pi_1 (E\setminus N,u) \ri B_p[D,K].$
This map is evidently a group homomorphism, and it is the braid
monodromy of $B$ with respect to  $E\tm D$, $\pi $, $u$.
We sometimes denote $\vp$ by $\vp_u.$\enddefinition

\definition{Definition} \ $\underline{\text{Braid monodromy of}\  B \
\text{with respect to} \ \pi, \, u.}$

When considering the braid induced from the previous motion as an element of
the group $B_p [\C_u, K]$ we get the homomorphism $\vp:\pi_1(E\setminus N,
u)\ri B_p[\C_u, K]$ which  is called the braid mondromy of $B$ with respect
to $\pi $, $u$.
We sometimes denote $\vp$ by $\vp_u.$\enddefinition
\mk
In order to present an example of a braid monodromy calculation, we recall  a
geometric model of the braid group and the definition of a half-twist.

\definition{Definition}\ $\underline{\text{Braid group}\ B_p[D,K]}$

Let $D$ be a closed disc in $\Bbb R^2,$ \ $K\subset D,$ $K$ finite.
Let $\Cal B$ be the group of all diffeomorphisms $\beta$ of $D$ such that
$\beta(K) = K\,,\, \beta |_{\partial D} = \text{Id}_{\partial D}$\,.
For $\beta_1 ,\beta_2\in \Cal B$\,, we
say that $\beta_1$ is equivalent to $\beta_2$ if $\beta_1$ and $\beta_2$ induce
the same automorphism of $\pi_1(D\setminus K,u)$\,. The quotient of $B$ by this
equivalence relation is called the braid group $B_p[D,K]$ ($p= \#K$).
We sometimes denote by $\overline\beta$ the braid represented by $\beta.$
The elements of $B_p[D,K]$ are called braids.
\enddefinition

\definition{Definition}\ \underbar{$H(\sigma)$, half-twist defined by
$\sigma$}

Let $D$, $K$ be as above. Let $a,b\in K\,,\, K_{a,b}=K\setminus \{a,b\}$ and
$\sigma$ be a simple (i.e. without self-intersections) path in $D\setminus \partial D$ connecting $a$ with $b$
such that $\sigma\cap K=\{a,b\}.$ Choose a small regular neighbourhood $U$ of
$\sigma$, $K_{a,b}\cap U=\emptyset$, and an orientation
preserving diffeomorphism $\psi :{\Bbb R}^2
\longrightarrow {\Bbb C}^1$ (${\Bbb C}^1$ is taken with usual ``complex''
orientation) such that $\psi(\sigma)=[-1,1]=\{z\in{\Bbb C}^1 \,|\, \text{Re}\
z\in [-1,1],\, \text{Im}z=0\, \}$ and $\psi(U)=\{z\in{\Bbb C}^1
\,|\,|z|<2\}$\,.
Let $\alpha(r),r\geqslant 0$\,, be a real smooth monotone function such that
$\alpha(r) = 1$ for $r\in [0,\tsize{3\over 2}]$ and $\alpha(r) = 0$ for
$ r\geqslant 2.$ Define a diffeomorphism $h:{\Bbb C}^1\to {\Bbb C}^1$ as
follows: for $z\in {\Bbb C}^1$, $z=re^{i\varphi }$ let $h(z)=re^{i(\varphi +
\alpha (r)\pi )}$. It is clear that the restriction of $h$ to
$\{ z\in {\Bbb C}^1\, \,
\mid \, \, |z|\leq {3\over 2}\, \} $ coincides with the positive rotation
on $\pi $, and that the restriction to $\{ z\in {\Bbb C}^1\, \,
\mid \, \, |z|\geq 2\, \} $ is the identity map. The diffeomorphism
$H(\sigma)=\psi^{-1}\circ h\circ \psi $ will be called a half-twist.
\enddefinition

A half-twist $H(\sigma)$ defines {\it a geometric braid}
$\overline \sigma$ (i.e. $p$ paths without self-intersections in
$D\times [0,1]$ starting at $K\times \{ 0\}$ and ending at
$K\times \{ 1\}$, $K=\{ q_1, \dots , q_p\}$). This braid can be presented as
$$ \alignat 2
&(\delta _j(t),t)=(q_i,t) \qquad \qquad \qquad && \text{if} \, \, q_j\neq
a,b; \\
&(\delta _j(t),t) =(\psi ^{-1}(e^{\pi it}),t) \qquad \qquad  && \text{if} \, \,
q_j= a; \\
&(\delta _j(t),t) =(\psi ^{-1}(-e^{\pi it}),t) \qquad \qquad  && \text{if} \,
\, q_j= b. \\
\endalignat
$$

The following is the basic braid monodromy associated to a single
singularity of an algebraic curve.

\proclaim{Proposition - Example  1.1} \ Let $E=\{x\in\C \, |\, |x|\leq
1\},$\ $D=\{y\in\C \, |\, y \leq R\},$ $R\gg 1, B$ is the curve
$y^2=x^\nu $,  $u=1$.
Clearly, here $n=2$, $N=\{0\} $, $K=\{-1, +1\}$ and $\pi_1(E-N, 1)$ is
generated by
$\G=\p E$ (positive orientation). Denote by $\vp:\pi_1(E\setminus N,
1)\ri B_2[D,K]$ the braid monodromy of $B$ with respect to $E\times D$,
$\pi $, $u$.

Then $\vp(\G)=H^\nu,$ where $H$ is the positive half-twist defined by $[-1,
1]$ (``positive generator" of $B_2[D, K]$).\endproclaim

\demo{Proof} \ We can write $\G = \{e^{2\pi it}, t\in[0,1]\}.$
Lifting $\G$ to $S$ we get two paths:
$$\align
\delta_1(t) &= \left(e^{2\pi it}, \ e^{2\pi i\nu t/2}\right)\\
\delta_2(t) &= \left(e^{2\pi it}, \ -e^{2\pi i\nu t/2}\right).\endalign$$

Projecting $\delta_1(t),$ $\delta_2(t)$ to $D$ we get two paths:
$$\alignat 2
&a_1(t) =e^{\pi it\cdot \nu}, \qquad \qquad  &&0\leq t\leq 1\\
&a_2(t) = -e^{\pi it\cdot \nu}, &&0\leq t\leq 1.
\endalignat$$

These paths induce a motion of \{1, -1\} in $D.$
This motion is the $\nu$-th power of the motion $\Cal M$:
$$\alignat 2
&b_1(t) =e^{\pi it}, \qquad \qquad  &&0\leq t\leq 1\\
&b_2(t) = -e^{\pi it}, &&0\leq t\leq 1.
\endalignat$$

The braid of $B_2\left[D, \{1, -1\}\right]$ induced by ${\M}$ coincides with
the half-twist $H$ corresponding to $[-1, 1]\subset D.$ Thus $\vp(\G)=H^\nu.$
\hfill$\qed$ \enddemo

We recall the notion of a geometric free base of the fundamental group of a
punctured disc and a basic property of the free base.

\definition{Definition}\ \underbar{A bush}.

Let $E$, $N=\{u_1,\dots , u_n\}$, $u$ be as above. Consider in $E$
ordered sets of simple paths $(T_1,\dots ,T_n)$ connecting $u_i$'s with $u$
such that
\roster
\item $T_i\cap T_j=u$ if $i\ne j$\,;
\item for a small circle $c(u)$ around $u$ each $T_i\cap c(u)$ is a single
point, namely $w_i$, and the order in $(w_1,\dots ,w_n)$ is consistent with
the positive (``counterclockwise'') orientation of $c(u)$\,.
\endroster
Let $c_j$ be the boundary of a closed disc $E_j$ of small radius with center
at $u_j$. Denote by $T_j'=T_j\setminus (T_j\cap E_j)$
and $l(T_j)=T_j'\cdot c_j \cdot {T'_j}^{-1}$ a loop (and the corresponding
element in $\pi _1(E\, \setminus N,u)$) in which a point moves
counterclockwise along
$c_j$.

We say that two such sets $(T_1,\dots ,T_n)$ and $(\widetilde T_1,\dots
,\widetilde T_n)$ are equivalent if
$$\ell(T_i)=\ell(\widetilde T_i)\quad\quad(\text{in }\,
\pi_1(E\setminus N,u)$$
for all $i=1,\dots ,n$.
An equivalence class of such sets is called a bush in $(E\setminus N,u)$\,.
The bush represented by $(T_1,\dots ,T_n)$ is denoted by
$\langle T_1,\dots ,T_n\rangle$\,.
\enddefinition

\definition{Definition}\ $\underline{\text{geometric base,}\ g-\text{base}}$

 Let $E,$ $N,$ $u$ be as above. A $g$--base of $\pi_1(E\setminus N,u)$ is
an ordered
free base of $\pi_1(E\setminus N,u)$ which has the form $(\ell(T_1),\dots
,\ell(T_n))$
where $\langle T_1,\dots ,T_n\rangle$ is a bush in $E\setminus N$\,.
\enddefinition

\midinsert
\medskip
 
\centerline{
\epsfysize=1.5in
\epsfbox{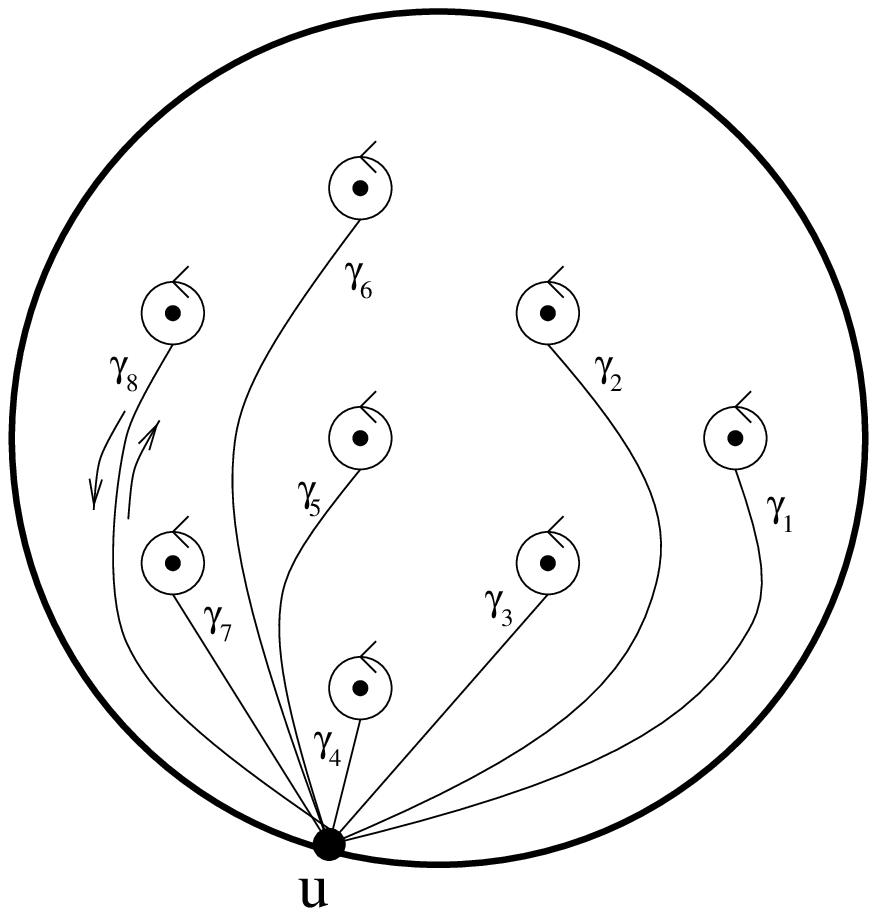}}  
\botcaption{Fig. 1.1}\endcaption
\endinsert

\proclaim{Proposition-Definition} Denote by $\wp$ the element of
$\pi_1(E\setminus N,u)$ represented by the loop $\partial E$ (with positive
orientation). There exists a unique element of $B_n[E,N]$\,, denoted by
$\Delta^2_n$ or $\Delta^2_n[E,N]$ such that for any $g$-base
$\Gamma_1,\dots,\Gamma_n$ of $\pi_1(E\setminus N,u)$ $$
(\Gamma_i)\Delta^2_n =\wp\Gamma_i \wp^{-1}\,.
$$
\endproclaim
\demo{Proof} \cite{MoTe1}, V.2.1.\edm
\remark{Remark}
Clearly, $\Dl_p^2$ acts as a full-twist around all the points of $N.$\endremark

\proclaim{Proposition 1.2} $\Dl_n^2\in\Center B_n[E,N].$\ep
\demo{Proof} \cite{MoTe1}, V.4.1.\edm
\proclaim{Proposition - Example 1.3} \ Let $B$ be a union of
$p$ lines, meeting in one point $s_0$, $s_0=(x(s_0),y(s_0))$.
Let $D$, $E$, $u$, $K=K(u)$ be as in the beginning of \S1. Let $\vp$ be
the braid monodromy of $B$
with respect to $ E\times D$, $\pi $, $u$. Clearly, here $N=$ single point
$x(s_0)$ and
$\pi_1(E\setminus N,u)$ is generated by $\G=\p E$.
Then $\vp(\G)=\Dl^2_p=\Dl^2_p \left[ D, K(u)\right].$\endproclaim

\demo{Proof} \ By a continuous change of $s_0$ and the $n$ lines passing
through $s_0$ (and by uniqueness of $\Dl^2_p$) we can reduce the proof to
the following case: $B=\cup L_k,$ \ $L_k\:$ $y=j_kx, \ j_k=e^{2\pi ik/p}$,
\quad $k=0,\dots ,p-1.$ Then $N=\{0\}.$ We can take $E=\{c| |x|\leq 1\},$
$u=1,$\ $\G=\p E=\{x=e^{2\pi it},$ $t\in[0,1]\}.$
Lifting $\p E$ to $B$ and then projecting it to $D,$ we get $p$ loops:
$$a_k(t)=e^{2\pi i(t+k/p)},\quad k=0,\dots ,p-1,\quad t\in[0,1].$$

Thus the motion of $a_k(0)$ represented by $a_k(t)$ is a full-twist which
defines the braid $\Dl^2_p\left[ D, \{a_k(0)\}\right] =
\Dl^2_p\left[ D, K(1)\right].$

(To check this last fact, see the corresponding actions in
$\pi_1\left(D\setminus K,u\right)$). \enddemo
\bigskip

\definition{Definition}\ \underbar{Frame of $B_p[D,K]$}.

Let $D$, $K=\{ q_1,\dots \, q_p\}$ be as in the beginning of \S1. Let us
choose a system of
simple smooth paths $(\sigma _1,\dots ,\sigma _{p-1})$ in $D\setminus
\partial D$
such that $\sigma _i$ connects $q_i$ with $q_{i+1}$ and
$L=\cup \sigma_j$ is a simple smooth path. The ordered system
of half-twists $(H_1,\dots , H_{p-1})$ defined by $\{ \sigma _i\}_{i=1}^{p-1}$
is called a frame of $B_p[D,K]$. Sometimes such a system of paths
$\{ \sigma _i\}_{i=1}^{p-1}$ will also be called a frame of $B_p[D,K]$.
\enddefinition

\proclaim{Theorem 1.4} Let  $(H_1,\dots ,H_{p-1})$ be a frame of $B_p[D,K]$.
Then $B_p[D,K]$ is generated by $H_1,\dots ,H_{p-1}$.
\endproclaim
\demo{Proof} See, for example, \cite{MoTe1}.\enddemo

\definition{Definition}\ \underbar{Generalized half-twist $\Dl _{i,j}$ given by
a skeleton $L_{i,j}=\bigcup _{r=i}^{j-1}\sigma _r$
.}

In the notation of the above definition let $\{ \sigma _r\}_{r=1}^{p-1}$ be a
frame of $B_p[D,K]$ and let $i<j$. Choose a small neighbourhood $U$ of the
path $L_{i,j}$ and an orientation preserving diffeomorphism
$\psi :{\Bbb R}^2 \longrightarrow {\Bbb C}^1$ such that
\roster
\item $\psi (L_{i,j})=[-1,1]=\{z\in{\Bbb C}^1 \,|\, \text{Re}\, z\in [-1,1],\,
\text{Im}\, z=0\, \}$ \,;
\item $\psi (U)=\{z\in{\Bbb C}^1 \,|\,|z|<2\}$\,;
\item the set $\{ \psi (q_r) \} _{r=i}^{j}$ is invariant relative to the
involution $\text{Re}\, z\mapsto -\text{Re}\, z$\,.
\endroster
Let $\alpha(r),r\geqslant 0$ be a real smooth monotone function such that
$\alpha(r) = 1$ for $r\in [0,\tsize{3\over 2}]$ and
$\alpha(r) = 0$ for $ r\geqslant 2.$
Define a diffeomorphism $h:{\Bbb C}^1\to {\Bbb C}^1$ as follows: for
$z\in {\Bbb C}^1$, $z=re^{i\varphi }$, let $h(z)=re^{i(\varphi +
\alpha (r)\pi )}$. The diffeomorphism $\Dl _{i,j}=\psi ^{-1}\circ h\circ \psi $
will be called a half-twist given by $L_{i,j}=\cup _{r=i}^{j-1}{\sigma _r}$.
\enddefinition
\remark{Remark} The full-twist $\Dl _{i,j}^2=(H_i\cdot \dots \cdot
H_{j-1})^{j-i+1}$. In particular, $\Dl _{i,i+1}=H_i$ and $\Dl _{1,p}^2 =
\Dl _{p}^2$, where the the half-twist $H_i$ and the
full-twist $\Dl _p^2=\Dl _p^2[D,K]$ was defined above.
\endremark

\subheading{\S 2.\ Positivity of the braid monodromy}

In this section, we shall show that the braid monodromy takes, in fact,
values in the
semigroup
$B_p^+$ of positive braid monodromy generated by the (counterclockwise)
half-twists.

We keep the same notations as in \S1, and introduce additional ones.

\definition{Definition} \ $\underline{\psi_T, \ \text{Lefschetz diffeomorphism
induced by a path} \ T }$

Let  $T$ be a path in $E\setminus N$ connecting $x_0$ with $x_1$,
$T:[0, 1]\ri E\setminus N$.
There exists a continuous family of diffeomorphisms $\psi_{(t)}:
D\ri D,\ t\in[0,1],$ such that $\psi_{(0)}=Id$, $\psi_{(t)}(K(x_0))=K(T(t)) $
for all $t\in[0,1]$, and  $\psi_{(t)}(y)= y$ for all $y\in\p D$.
For emphasis we write $\psi_{(t)}:(D,K(x_0))\ri(D,K(T(t))$.
Lefschetz diffeomorphism induced by a path $T$ is the diffeomorphism
$$\psi_T= \psi_{(1)}: (D,K(x_0))\usr (D,K(x_1)).$$\enddefinition
\mk
Since $ \psi_{(t)} \left( K(x_{0})\right) = K(T(t))$ for all $t\in
[0,1]$, we have a family of
canonical isomorphisms
$$\psi_{(t)}^\nu: B_p\left[ D, K(x_{0})\right] \usr B_p\left[
D, K({T(t)})\right], \ \quad \text{for all} \, \, t\in[0,1].$$

\definition{Definition} \ $\underline{L_T, \ \text{Lefschetz isomorphism
induced by} \ T}$
$$\align
&L_T= \psi_T^\nu=\psi_{(1)}^\nu: B_p\left[  D, K(x_{0})\right] \usr
B_p  \left[  D, K(x_1)\right]\\
&(\overline{\be}) L_T = \overline{\psi\1_T\circ \be \circ \psi_T}.\endalign$$

It is easy to check that $L_T$ depends only on the homotopy class of
$T.$\enddefinition

\demo{Notation}\ $\underline{\psi_{T, B}\ \ L_{T,B}}$

$L_T$ and $\psi_T$ depend not only on $T,$ but also on $B.$
To avoid confusion, we shall emphasize and use
$\psi_{T, B}$ and $L_{T,B}.$\edm

\remark{Remark} \ There is an equivalent definition of braid monodromy
$\vp_u:\pi_1(E\setminus N, u)\ri B_p [ D, K].$
Take any $\delta\in\pi_1(E\setminus N, u).$
Let $\underline\delta$ be a loop representing $\delta.$
Then $\vp(\delta)$ is the braid represented by the diffeomorphism
$\psi_{\underline\delta}.$\endremark

\demo{Proof}\ Both $\psi_{\underline\delta}$ and $\vp (\delta)$ are induced
from ``pushing'' along $\G.$\enddemo

\proclaim{Lemma 2.1}
\roster
\item $\psi_{T_1 T_2} = \psi_{T_1} \circ \psi_{T_2}.$
\item $H(\s) L_t = H((\s)\psi_T).$
\item If $T$ is a closed loop, then $\psi_T$ defines a braid
$\overline{\psi}_T$ in $B_p$ and $(b) L_T = \overline{\psi}_T\1 b\ov{\psi}_T$
(composition from left to right).\endroster\ep

\demo{Proof}\
 (1) and (2)\ are proved immediately.

(3)\ Assume $b=\overline{\be}.$
By definition of $L_T,\ (\overline{\be}) L_T =
\overline{\psi_T\1\circ\be\circ\psi_T}.$
Since $\psi_T$ and $\be$ each defines a braid, $(b) L_T = \ov{\psi}_T\1
\circ\overline{\be}\circ \ov{\psi}_T =
\overline{\psi}_T\1 b\overline{\psi}_T.$\hfill$\qed$  \edm

Let $M'=\{s_j\}_{j=1}^n$ be the singular points of $\pi\bigm|_B$,
$\pi(M')=N.$ For every $j=1,\dots , n$,
let $D_j'$ be  a small disk on $y$-axis centered at $y(s_j)$ such that
$D_j\subseteq D$ and that $(x(s_j)\times D_j' \cap B=s_j.$
Let $E_j'$ be a sufficiently small closed neighborhood of $x_j= x(s_j)$ on the
$x$-axis such that $E_j'\cap N=x_j$\ and
the number $\#(x\times\Int (D_j')\cap B)$ is independent of $x$ for all
$x\in E_j'\setminus \{ x_j\} $.
We call this number the local degree of $\pi$ at $s_j$ or $\deg_{s_{j}} \pi.$
Let $m_j= \deg_{s_{j}}  \pi.$ Choose a point $x_j'\in\partial E.$
Let $K'(x_j')=K(x_j')\cap D_j'.$

\definition{Definition} \
$\underline{\psi_{T,s_j}, \ \text{Lefschetz embedding induced by}\ s_j\
\text{and} \ T}$

  Let $T:[0,1]\ri\C$ be a path in
$E\setminus (N\cup (\text{Int} \ E'_j))$ connecting $x_j'$ to a point $u'\in
E\setminus N$. The diffeomorphism
  $\psi_{T,s_{j}}=\psi_T\bigm|_{ D_{j}'}:(D_j',K'(x_j'))\ri(D,K(u'))$ is
called Lefschetz embedding induced by  $s_j$ and $T$, where
  $\psi_T$ is the Lefschetz diffeomorphism induced by $T$. \enddefinition

\remark{Remark} Let $m_j=\deg_{s_{j}}\pi.$ Take $m_j$ liftings of $T$ to $B$
starting at the different points of $K'(x_{j}').$
These liftings are real curves in $T\times D.$
We can think of $\psi_T$ as ``pullings" of $K'(x_{j}')$ in $T\times D$ along
these real curves.\endremark

\definition{Definition} \ $\underline{ L_{T,s_{j}}, \ \text{Lefschetz
injection induced by} \ T}$

Let $s_0\in N.$ Let $D_0'$, $x_0$, $E_0$, $m_0$, $u^{\prime }$ be as in the
definition of
Lefschetz embedding. Consider $\psi_{T,s_{j}}:(D_j'\times K(x_j'))
\ri(D,K(u'))$ a Lefschetz embedding induced by $s_j$ and $T.$
We have $\psi_{T,s_{j}}(K'(x'_j))\subset K(u')$ and
$(K(u')\setminus \psi_{T,s_{j}}(K(x_j'))\cap\psi_{T,s_{j}} (\Int  D_j')=\phi.$
The Lefschetz injection induced by $T$ is the canonical injection induced from
$\psi_{T,s_{j}},$\ $L_T=L_{T,s_{j}}=\psi^\nu_T: B_{m_{j}}\left[D_j',
K'(x_{j}')\right]\hookrightarrow B_p[D,K(u')]$ which is
is well-defined by the above inclusions.\enddefinition

To compute the braid mondromy, we need to know
$\{\vp(\delta_j)\}^n_{j=1}$ for the
$g$-base $\{\delta_j=\ell(\gamma_j)\}_{j=1}^n$ defined by the bush
$\{\gamma_j\}_{j=1}^n$ in $(E\setminus N,u).$ We actually
need to be able to compute $\vp(\delta_j)$ for such $\{\delta_j\}$ since it is
the basic data for any applications of braid monodromy. We can represent each
$\delta_j$ as $\tg_j^{-1}\circ\p E'_j\circ\tg_j$ ($E'_j$ was defined earlier).
Thus, to know $\vp(\delta_j)$ it is sufficient to know the Lefschetz injection
$L_{\tg_{j}}: B_{m_{j}} [ D'_j, K'(x'_j)]\ri B_p[ D, K]$ and also the
local braid monodromy $\vp_{s_{j}}$ of $B \cap (E'_j\times D'_j)$ with
respect to
$E'_j\times D'_j,$ $\pi_, x'_j$ as defined here.

\definition{Definition} \
$\underline{\vp_{s_{j}} \ \text{local braid monodromy of} \ B \ \text{at} \
s_j}$

The local braid monodromy of $B$ at $s_j$ is
$$\vp_{s_{j}}:\pi_1(E'_j\setminus \{ x_j\} ,x'_j)\ri B_{m_{j}}\left[ D'_j,
K'(x'_{j})\right],$$
the braid monodromy of
$B\cap[E_j'\times D_j']$ with respect to $E'_j\times D'_j$, $\pi $, $x'_j$
 ($K'(x'_{j})=K(x'_j)\cap D_j'$).
\mk
It is clear that $\vp_{s_{j}}$ is determined only by $\vp_{s_{j}}(\p E'_j).$
The following lemma is evident.\enddefinition

\proclaim{Lemma  2.2} \ Let $T$ be any path
in $E\setminus (N\cup (\Int \ E'_j))$ connecting $x'_j$ with $u$ and
$\dl=T^{-1}\circ\p
E'_j\circ T=\ell(T).$ Then $\vp(\delta) = \left(\vp_{s_{j}} (\p E'_j)\right)
L_T.$ In particular, $\vp(\delta_j) = \left(\vp_{s_{j}} (\p E'_j)\right)
L_{\tg_{j}}.$\endproclaim

\remark{Remark} The   Lemma actually indicates that the
braid monodromy $\vp$ is completely determined if we know local
braid monodromies and Lefschetz injections for some bush in $E\setminus
N.$\endremark

\proclaim{Lemma 2.3} \ Let $s_j$ be a singularity of $B$ which is
locally presented by $y^2=x^\nu,$ that is, $m_j= \deg_{s_{j}} \pi=2$.
Then
$\vp(\delta _j)=(H_{j})^\nu$, where $H_j$ is a positive half-twist defined
by some path $\sigma $,  and in particular it is a positive braid.
\endproclaim

\demo{Proof} \ Follows from Proposition - Example  1.1. and Lemma 2.2.\enddemo

\proclaim{Proposition 2.4} \  Let $\vp:\pi_1(E\setminus N,u)\ri B_p[D,K]$
be  braid monodromy of a curve $B$, $\{\delta_j\}$ a $g$-base of
$\pi_1(E\setminus N,u).$
Then all $\vp(\delta_j)\in B^+_p= B_p^+[D,K].$\endproclaim

\demo{Proof} \ Given a curve $B$, we can find a curve $B^{(1)}$ close enough
to $B,$ nonsingular and of the same degree.
Let $K^{(1)} = \{y\bigm| (u,y)\in B^{(1)}\}$ and $M^{(1)},$\ $N^{(1)}$ be
as in \S1. We can naturally identify $B_p[  D, K]$ with $B_p[  D,K^{(1)}].$
Each $s_j\in M$ splits into a finite set of singular points
$\{s_{j_i}\}\subseteq
M^{(1)}$ which locally are of type $x=y^2.$
Each $x_j\in N$ will split into points $\{x_{j_i}\} = \{\pi(s_{j_i})\}\subseteq
N^{(1)}.$ Clearly, $N^{(1)} =\{x_{ji}\}_{j,i}.$
Let $\vp^{(1)}$ be the braid monodromy of $B^{(1)}$ with respect to
$E\times D$, $\pi $, $u$, $N^{(1)}$. We can find $\{\delta_{j_i}\}$ a
$g$-base of
$\pi_1(E\setminus N^{(1)}, u)$ such that each $\delta_j=\prod\limits_i
\delta_{j_i}.$\ Natural identification of $B_p[D,K]$ and
$B_p[  D, K^{(1)}]$ will give us that each $\vp(\delta_j)=\prod\limits_i
\vp^{(1)}(\delta_{j_i})$ (we use the fact that $B^{(1)}$ is very close to $B$).
By Lemma 2.3, each $\vp^{(1)} (\delta_{j_i})$ is a positive half-twist.
Thus, each $\vp(\delta_j) = \prod\limits_i \vp^{(1)} (\delta_{j_i}) \in B^+_p.$
\hfill $ \qed$\edm
\bigskip
\subheading{\S 3.\ Braid monodromy  of a projective curve}

\definition{Definition} \ \underbar{Braid monodromy of a projective curve}

Let $B$ be an algebraic curve of degree $p$ in $\C \Bbb P^2.$
Choose generically a line $L$ at infinity  $(\#(L\cap B) = p)$ and affine
coordinates $(x,y)$ in $\C^2 = \C \Bbb P^2\setminus L$, so that the
projection $\pi (x,y)=x$ on $x$-axis of the curve
$B\cap \C^2$ is generic (in particular, the center of this projection
in $\C \Bbb P^2$ must be outside of $B$). Let $\pi(x,y)=x.$
Let $N=\{x\in\C \, | \, \pi^{-1}(x)\cap B\lvertneqq p\},$ $E$  is a closed
disk on
the $x$-axis with $N\subset\Int (E) $, $D$ is a closed disk on the $y$-axis
with $\pi^{-1} (E)\cap B\subset E\times D.$ Choose $u\in\p E.$

The  braid monodromy of $B$ with respect to $L$, $u$ is the braid monodromy of
$B\cap (E\times D)$ with respect to  $E\times D$, $\pi $, $u$, i.e., the
homomorphism
$$\vp:\pi_1(E\setminus N, u)\ri B_p [ D, K]$$
defined in \S1.\enddefinition

\proclaim{Proposition 3.1} \ Let $B$ be an algebraic curve of degree $p$ in
$\CP^2$. Let $L$, $\pi $, $u$, $D$, $E$, $K(u)$ be as above.
Let $\vp$ be the braid monodromy of $B$ with respect to  $L$, $\pi $, $u$. Let
$\delta_1,\dots,\delta_q$ be a $g$-base of $\pi_1(E\setminus N,u).$ Then
$$\prod^q_{i=1} \vp(\delta_i)=\Dl^2_p=\Dl^2_p [u\times D, K(u)].$$
\endproclaim

\demo{Proof} \ Because $\prod\limits^q_{i=1} \ \delta_i=\p E$ (positive
oriented), we have to prove that $\vp(\p E)=\Dl^2_p.$
We can assume $E$ arbitrarily big, so that $\p E$ will be very close to
$\infty$ at the $x$-axis.
Continuously deforming coefficients of the equations of $B$ such that new
curves will be transversal to $L_{\infty }$ we can reduce the proof to an
equation which
defines union of $n$ lines intersecting at a single point.
Now use Proposition - Example 1.3. \hfill $\qed$
\enddemo

\proclaim{Lemma 3.2} $\Dl_p^2\in B_p^+.$\ep

\demo{Proof} Proposition 2.4.\edm

\definition{Definition} \
$\underline{\text{Braid monodromy factorization}}$ (associated to
projective curve)

Braid monodromy factorization associated to a plane  projective curve is
a product  of the form $\Dl^2_p=\prod\limits_i\vp (\delta_i),$ where
$\vp $ is the braid monodromy of the projective curve and $\{\dl_j\}$ is
a $g$-base of $\pi_1(E\setminus N_1,u).$\enddefinition

\remark{Remarks}

(1) \ A $g$-base of $\pi_1(E\setminus N,u)$ and the corresponding product-form
determine the braid monodromy. For applications it is usually sufficient to
know a product-form of a braid monodromy without reference
to a particular $g$-base.

(2) \ The product-form is not a prime factorization of $\Dl^2_p$ unless $B$
is nonsingular. For a nonsingular $B$, each $\vp(\dl_i)$ is a
positive half-twist which is a prime element of $B^+_p.$\endremark

\proclaim{Proposition 3.3} \ Let $B$ be a (generalized) cuspidal curve on
$\C \Bbb P^2$ (that is, all singularities of $B$ are locally given by
$y^2=x^{\nu }$, $\nu \in \Bbb N$). Then any product-form of the braid
monodromy of $\delta$ can be written as $\Dl^2_p=\prod\limits_i
(Q_i^{-1} H^{\nu_i}_1 Q_i)$ where $H_1$ is a
positive half-twist and each $\nu_i\in \Bbb N.$\endproclaim

\demo{Proof} \ Recall that we are using generic projections of
$\C^2\overset\pi\to\rightarrow \C$ with respect to  the projective curves. Each
singularity of $\pi|_B$ is of the type $y^2=x^\nu,$ $\nu \in \Bbb N$. Now use
Lemma 2.3 to get $\vp(\dl_i) = H_i^{\nu_i}$, where $H_i$ is  a half-twist.
Every two half-twists in $B_p$ are conjugate, so for all $i$ there exist $Q_i$
such that
$H_i=Q_i\1 H_1 Q_i.$
Thus, $\Dl^2_p=\Pi \vp(\dl_i) = \Pi Q_i\1 H_1^{\nu_i} Q_i.$

\hfill $\qed$
\enddemo

\remark{Remark}  We can take any half-twist for $H_1$.
\endremark

\subheading{\S 4.\ Braid monodromy factorizations of $\Dl_p^2$ and Hurwitz
equivalence}

 From Proposition 3.1 and  2.4, we know that a braid monodromy factorization
$$\Dl^2_p[D,K] =\prod\limits_i\vp(\delta_i)$$ is a factorization of $\Dl^2_p
[D,K]$ in $B^+_p[D,K]$ induced from a $g$-base of $\pi_1(E\setminus N,u)$.

We define an equivalence relation on the set of $B_p^+$-factorizations of
$\Dl_p^2.$ We start by classifying a Hurwitz move on $G\times\dots\times G$
($G$ is a group) or on a set of factorizations.

\definition{Definition} \ $\underline{\text{Hurwitz moves}\ R_k,\ R^{-1}_k\
 \text{on}\ G^m}$:\quad

Let $\vec t= (t_1,\dots ,t_m)\in G^m$\,. We say that
$\vec s =(s_1,\dots ,s_m)\in G^m$ is obtained from
$\vec t$ by the Hurwitz move $R_k$ (or $\vec t$ is
obtained from $\vec s$ by the Hurwitz move $R^{-1}_k$) if
$$
\align
&s_i = t_i \quad\text{for}\  i\ne k\,,\, k+1\,,\\
&s_k = t_kt_{k+1}t^{-1}_k\,,\\
&s_{k+1} =t_k\,.
\endalign
$$\enddefinition
\definition{Definition}\ \underbar{Hurwitz move on a factorization}

Let $G$ be a group $t\in G.$  Let  $t=t_1\cdot\dots\cdot t_m=
s_1\cdot\dots\cdot s_m$ be two factorized
expressions of $t.$ We say that $s_1\cdot\dots\cdot s_m$ is obtained from
$t_1\cdot\dots\cdot t_m$ by a Hurwitz move $R_k$ if $(s_1,\dots ,s_m)$ is
obtained from $(t_1,\dots ,t_m)$ by Hurwitz move $R_k$\,. \enddefinition

\definition{Definition}\ \underbar{Hurwitz equivalence of factorization}

Two factorizations are Hurwitz equivalent if they are obtained  from each
other by a
finite sequence of Hurwitz moves. \enddefinition

In order to study equivalence relations, we first state that braid
monodromy and
Hurwitz moves are commutative.
\proclaim{Lemma 4.1} \rom{(Proof in \cite{MoTe1}, Chapter II)}.
 Let $D,K,u$ be as above.
\roster
\item"{(a)}" if $\Gamma_1,\dots ,\Gamma_n$ is a $g$--base of
$\pi_1(D\setminus K,u)$\, then $\Gamma_1\cdot\dots\cdot\Gamma_n$ is
represented by the
loop $\partial D$ (taken with positive orientation).
\item"{(b)}" If $\{\Gamma^{\prime}_i\}$ and $\{\Gamma_i\}$ are two
$g$--bases of $\pi_1(D\setminus K,u)$ then each $\Gamma^{\prime}_i$ is
conjugate
to some $\Gamma_{j_i}$\,.
\item"{(c)}" By applying a Hurwitz move to a $g$--base, we get a
$g$--base (see \text{\rm{Fig. 4.1)}}.
\item"(d)" Any two $g$--bases can be obtained  from each other by a finite
sequence of Hurwitz moves.\endroster

\midinsert
\medskip
 
\centerline{
\epsfysize=1.2in
\epsfbox{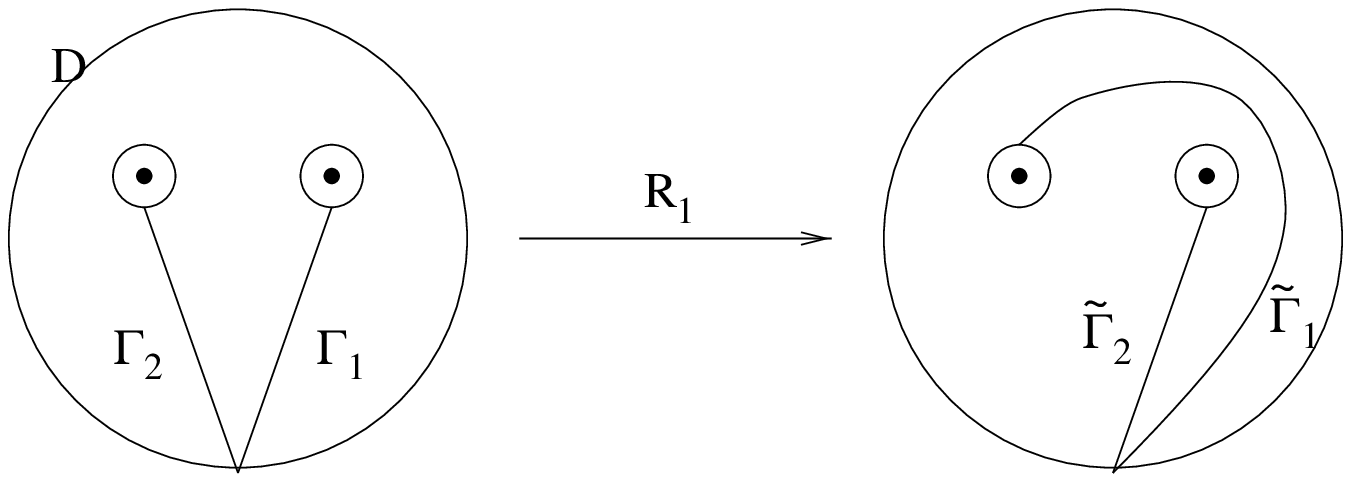}}  
 
\botcaption{Fig. 4.1}\endcaption
\endinsert

 \endproclaim

\proclaim{Proposition 4.2} \ Let $B$ be  a curve in $\C^2$, and let
$N$, $E$, $D$, $\pi $, $u$ be as above. Let $\varphi$ be the braid
monodromy of $B$
with respect to
$E\times D$, $\pi $, $u$. Let $n=\#N$. The following diagram
$$\CD
\{g\text{-bases of} \
\pi_1(E\setminus N,u)\} @>\vp^n>> (B_p[D,K])^n\\
@VR_k VV              @VR_k VV\\
\{g\text{-bases of} \
\pi_1(E\setminus N,u)\} @>\vp^n>>
(B_p[D,K])^n\endCD$$
\flushpar where $R_k,$ is the $k$-th Hurwitz move, is commutative.
\endproclaim

\demo{Proof} \ Follows immediately from definitions and the fact that
$$\vp\left(\delta_k\delta_{k+1}\delta_k^{-1}\right)
=\vp(\delta_k)\vp(\delta_{k+1}) (\vp(\delta_k))^{-1}. \hfill \qed$$\enddemo

Different braid monodromy factorizations related to $B$ that are derived from
different $g$-bases are equivalent to each other:
\proclaim{Lemma 4.3} \ Let $B$ be a curve, let $L$, $\pi $, $E$, $D$, $N$,
$u$, $K$ be as
in \S3. Let $\vp$ be the braid monodromy of $B$ with respect to $L$, $u$.
Let $\Dl^2_p[D,K]=\prod\limits_i\vp(\delta_i)$ and $\Dl^2_p[D,K]=
\prod\limits_i\vp(\delta'_i)$ be two braid monodromy factorizations
of $\Dl^2_p[D,K]$ corresponding to $\vp$ and two $g$-bases of
$\pi_1(E\setminus N,u).$ Then the two factorizations are Hurwitz
equivalent.\endproclaim

\demo{Proof} \ Two $g$-bases $\{\delta_i\}$ and
$\{\delta'_i\}$ of $\pi_1(E\setminus N,u)$ can be obtained
  from each other by a finite sequence of Hurwitz moves (Lemma 4.1).
By Proposition 4.2 the same sequence of Hurwitz moves will transform
$\{\vp(\delta_i)\}$ into $\{\vp(\delta'_i)\}.$
Thus the factorizations $\prod\limits_i\vp(\delta_i)$ and $
\prod\limits_i\vp(\delta'_i)$ are equivalent. \hfill $\qed$\enddemo

\proclaim{Lemma 4.4} \ In the notation of the previous Lemma, if $\prod_i
Z_i$ is Hurwitz equivalent to $\prod \vp(\delta_i)$, then there exists a
$g$-base $\{\delta_i'\}$ of $\pi_1(E\setminus N,u)$ such that
$Z_i=\vp(\delta'_i).$\endproclaim

\demo{Proof} \ Let $\ve$ be the sequence of Hurwitz moves that takes
$\vp (\delta_i)$ to $\{Z_i\}.$ Apply $\ve$ on $\{\delta_i\}$ to get
$\{\delta'_i\}.$ By Proposition 4.2,  $\vp(\dl'_i) = Z_i.$
\hfill $ \qed$\enddemo
 We conclude that:
\proclaim{Theorem 4.5} \ Let $B$ be a projective curve in $\Bbb C\Bbb P^2$.
Let $\vp\:\pi_1(E\setminus N,u)\ri B_p[D,K]$ be its braid monodromy.
The set of all braid monodromy factorizations of $\Dl^2_p[D,K]$ associated
to $B$ (presented by $\Dl^2_p=\prod\limits_i\vp(\delta_i),$  where
$\{\delta_i\}$ is a $g$-bases of $\pi_1(E\setminus N,u))$ occupy a full
equivalence class of factorizations of $\Dl^2_p$ in $B^+_p$.\endproclaim

Let $B$ be an algebraic curve in $\Bbb C^2$ and let
$\vp\:\pi_1(E\, \setminus \, N,u)\ri B_p[D,K]$ be the braid monodromy of $B$
defined by a braid monodromy factorization
$\Dl^2_p=\prod\limits_i\vp(\delta_i)$, where $\{\delta_i\}$ is a
$g$-base of $\pi_1(E\setminus N,u))$. Acting on $(D,K)$ by a diffeomorphism
$\beta $, we obtain a new braid monodromy factorization
$\Dl^2_p=\prod\limits_i\beta ^{-1} \vp(\delta_i)\beta$ associated to $B$.
\definition{Definition}\ \underbar{Braid monodromy factorization type}.

Two braid monodromy factorizations are called {\it equivalent with respect to
Hurwitz moves and conjugations} if one of them can be obtained from the other
by a finite sequence of Hurwitz moves, followed by a simultaneous
conjugation by an element $\beta \in B_{p}$.

Two braid factorizations belong to the same {\it braid monodromy
factorization type} if they are Hurwitz and conjugation equivalent.
\enddefinition

\subheading{\S 5.\ The centralizer of the multiple half-twists}

In this section we give a description of generators of the centralizer of the multiple
half-twist which we will use in the proof of Theorem 1.
We keep the same notations as in \S1.

\proclaim{Theorem 5.1} Let $(H_1,\dots ,
H_{p-1})$ be a frame of $B_p[D,K]$ given by a system of paths
$\{ \sigma _i\}_{i=1}^{p-1}$ and let $X=H_1^{\nu }$, $\nu \in \Bbb N$. Then
the centralizer $C(X)$ of $X$ in $B_p[D,K]$ is generated by
 $\Dl _{1,j}^2$, $j=3,\dots , p,$ and $H_j$, $j=1,3,\dots , p-1$.
\endproclaim
\demo{Proof}
Denote by $C_1(X)$ a subgroup of $B_p[D,K]$ generated by
$H_1$, $\Dl _{1,j}^2$, $j=3,\dots , p$, and the set of all half-twists
$H(\sigma )$ given by simple paths $\sigma $ starting and ending at
$K$ and non-intersecting with $\sigma _1$. Let $C_2(X)$ be a subgroup
generated by $\Dl _{1,j}^2$, $j=3,\dots , p$, and $H_j$, $j=1,3,\dots , p-1$.
It is clear that $C_2(X)\subset C_1(X)\subset C(X)$ and we must prove the
inverse inclusions.

It is sufficient to prove Theorem 5.1 for even $\nu $, since $C(X)\subset
C(X^2)$.

Without loss of generality,  we can assume that
$$D=\{ \, (v_1,v_2)\in \Bbb R^2\,
|\, v_1^2+v_2^2\leq p^2\, \} , \qquad K=\{ q_0=(0,0),\dots ,\, q_{p-1}=
(p-1,0)\},$$
$$\sigma _i=[i-1,i]=\{ \, (v_1,v_2)\in \Bbb R^2\,
|\, i-1\leq v_1\leq i,\, \, v_2=0\, \} .
$$
Denote by $\sigma _0=[-p,0]=\{ \,
(v_1,v_2)\in \Bbb R^2\,
|\, -p\leq v_1\leq 0,\, \, v_2=0\, \}$ and $\sigma _p=[p-1,p]=
\{ \, (v_1,v_2)\in \Bbb R^2\, |\, p-1\leq v_1\leq p,\, \, v_2=0\, \}$.
We choose a point $u_0\in \partial D$ such that $u_0$ does not lie in the lines
$\{ v_1=i\}$, $i=0,\dots , p-1$, and $\{ v_2=0\}$.
Consider an element $\gamma \in \pi _1(D\setminus K, u_0)$, $\gamma :[0,1]
\to D\setminus K$, \, $\gamma (0)=\gamma (1)=u_0$. By slightly changing
$\gamma$, we can assume that it is in a general position with respect to the
lines  $\{ v_1=i\}$, $i=0,\dots , p-1$, Ï $\{ v_2=0\}$. The coordinates
$(v_1,v_2)$ define an orientation on $D$ and on these lines. For each point
$y\in \gamma \cap L$, where $L$ is an oriented line in general position
with respect to
$\gamma $, the orientation chosen above allows us to define an intersection
index
$(\gamma ,L)_y$ equal to $\pm 1$.
Let $t_0=0<t_1<\dots <t_{n-1}<t_n=1$
be a sequence of $t\in [0,1]$ for which $\gamma (t)$
belongs to one of the lines considered above. We associate a sequence
$c(\gamma )=(a_0,\dots , a_n)$ of symbols $a_i\in \{u_0^{\pm 1}, o_0^{\pm
1},\dots ,o_p^{\pm 1}, h_0^{\pm 1},
\dots ,
h_{p-1}^{\pm 1},l_0^{\pm 1},\dots ,l_{p-1}^{\pm 1}\} $ to this loop (a {\it
code} of $\gamma$) as follows
\roster
\item $a_0=u_0$ and $a_n=u_0^{-1}$ \, ;
\item $a_i=o_j^{\pm 1}$ if $\gamma (t_i)\in \sigma _j$ and the power
coincides with $(\gamma ,\sigma _j)_{\gamma (t_i)}$ \,;
\item $a_i=h_j^{\pm 1}$ ifÏ $\gamma (t_i)\in L_j$, where
$L_j=\{ v_1=j, v_2>0\} $, and the power coincides withÏ
 $(\gamma ,L_j)_{\gamma (t_i)}$ \,;
\item $a_i=l_j^{\pm 1}$ if $\gamma (t_i)\in L_j$, where
$L_j=\{ v_1=j, v_2<0\} $, and the power coincides withÏ
 $(\gamma ,L_j)_{\gamma (t_i)}$.
\endroster

The code $c(\gamma )=(a_0,\dots , a_n)$ is said to be {\it reduced} if
 $a_i\neq a_{i+1}^{-1}$ for each $i$.
We can associate the reduced code $c_r(\gamma )$ to a code $c(\gamma)=(a_0,
\dots , a_n)$ removing from $(a_0,\dots , a_n)$
all pairs $a_i,a_{i+1}$ such that $a_i=a_{i+1}^{-1}$.

We call $l_c(\gamma )=|c(\gamma )|=n+1$ a {\it length} of $c(\gamma
)=(a_0,\dots , a_n)$.

Similarly, to each path $\sigma $ connecting points $q_i$ and
$q_j$ in $D\setminus (K\setminus \{ q_i,q_j\} )$, we can associate a code
$c(\sigma )$ if we add symbols $q_l^{\pm }$, $l=1,\dots p$,  to the symbols
defined above, and we can define a notion of reduced code satisfying the
following conditions:
\roster
\item $a_1=q_i$ and $a_n=q_j^{-1}$\, ;
\item $a_2\in \{ h_{i-1}^{-1},l_{i-1}^{-1},h_{i+1},l_{i+1} \}$  \, ;
\item $a_{n-1}\in \{ h_{j-1},l_{j-1},h_{j+1}^{-1},
l_{j+1}^{-1} \}$  \, ;
\item $a_l\neq a_{l+1}^{-1}$ for each pair $a_l,a_{l+1}$ \, .
\endroster

We define a sign $sgn(\sigma )$ of a code $c_r(\sigma)$ putting
$sgn(\sigma )=1$ if $a_2\in \{ h_{i-1}^{-1},l_{i-1}^{-1}\}$ and $sgn(\sigma
)=-1$ if $a_2\in \{ h_{i+1},l_{i+1} \}$.
\proclaim{Lemma 5.2} Loops $\gamma _0$ and $\gamma _1$ (resp. paths $\sigma
^{\prime }$ and $\sigma ^{\prime \prime }$) are homotopic in
$D\setminus K$ if and only if
$c_r(\gamma _1)=c_r(\gamma _2)$ (resp. $c_r(\sigma ^{\prime })=
c_r(\sigma ^{\prime \prime })$).
\endproclaim
\demo{Proof}
To each code $c(\gamma )=(a_0,\dots , a_n)$
we associate the element $a_0\cdot \dots \cdot a_n$ of a free group $F$
generated by
$$u_0, o_0,\dots ,o_p, h_0,\dots ,
h_{p-1},l_0,\dots ,l_{p-1}, q_1,\dots , q_p .$$
It is clear that a reduced word in $F$ corresponds to a reduced code.
It is well known that for each element in $F$ the reduced word representing
this element is uniquely defined. Therefore the reduced code is uniquely
defined for each code.

It is clear that if $\gamma _0$ and $\gamma _1$ are homotopic in
$D\setminus K$, then there exists a homotopy $\gamma _s$, $s\in [0,1]$, such
that for almost all $s$ except a finite number of
$s\in \{ s_1,\dots , s_k\}$, the loops
$\gamma _s$ are in general position with respect to the lines
$\{ v_1=i\}$, $i=0,\dots , p-1$, Ï $\{ v_2=0\}$, and for
$s\in \{ s_1,\dots , s_k\}$ the loops
$\gamma _s$ are touching one of these lines at one of the intersection
points,
and meet transversally at the other intersection points. Hence  homotopic
loops
have the same reduced code. The inverse statement that if the reduced codes
of homotopic loops $\gamma _0$ and $\gamma _1$ are equal to each other, then
$\gamma _0$ and $\gamma _1$ are homotopic is evident.

The case of two paths is similar to the case of two loops considered above.
\enddemo

First, we will show that $C_1(X)=C(X)$. Consider an element
$G\in C(X)$.
\proclaim{Lemma 5.3} Let a diffeomorphism $h$ be a representative of
$G\in C(X)$. Then the simple paths $\sigma =h(\sigma _1)$ and $\sigma _1$
(considered as non-oriented paths) are homotopic in $(D,K)$ (a homotopy
leaving fixed $K$).
\endproclaim
\demo{Proof} We have
$$X=G^{-1}XG=(G^{-1}H(\sigma _1)G)^{\nu }=
H(h(\sigma _1))^{\nu }=H(\sigma )^{\nu },$$
i.e. $X=H_1^{\nu }$ can be represented as $\nu$-th power of the half-twist
defined by  $\sigma $.
Therefore Lemma 5.3 follows from

\proclaim{Lemma 5.4} Let multiple full-twists
$H(\sigma _1)^{\nu}$ and
$H(\sigma )^{\nu}$ represent the same element in $B_p[D,K]$,
where $\sigma $ is a simple path. Then $\sigma $ and $\sigma _1$
(considered as non-oriented paths) are homotopic in $(D,K)$ (a homotopy
leaving fixed  $K$).
\endproclaim
\demo{Proof} Consider the reduced code
$c_r(\sigma )=(q_i,a_2, \dots,a_{n-1}, q_j^{-1})$. If the length
$|c_r(\sigma )|=2$, then  $\sigma $ is homotopic to a path $\sigma _l$
belonging to the chosen frame. Therefore, $l=1$ follows from
$H(\sigma _1)^{\nu }=H(\sigma _l)^{\nu}$.

Let us show that $|c_r(\sigma )|<3$. In fact, assume that
$|c_r(\sigma )|\geq 3$. Then for some $s\neq 0,1$, there is a symbol
$a_{l_0}\in \{
h_s^{\pm 1}, l_s^{\pm 1}\}$   in the reduced code  $c_r(\sigma )$.
Let we have $a_{l_0}=h_s^{\varepsilon }$, $\varepsilon =\pm 1$
(the case $a_{l_0}=l_s^{\varepsilon }$ is similar and it will be omitted). We
choose a point $u_0$ such that $s-1<v_1(u_0)<s,\, \, v_2(u_0)>0$
and consider a loop $\gamma \in \pi _1(D\setminus K, u_0)$ whose reduced
code  $c_r(\gamma )=(u_0, o_{s-1},l_s, o_s^{-1},h_s^{-1},
u_0^{-1})$ ($\gamma $ coincides with the $(p-s)$th element of a  $g$-base).
Then $(\gamma )H(\sigma _1)^{\nu}=\gamma $ and, therefore,
$$c_r((\gamma )H(\sigma _1)^{\nu })=c_r(\gamma ),$$
since we can choose a loop representing  $\gamma $ which does not intersect
$\sigma _1$.
To find a code $c((\gamma )H(\sigma )^{\nu })$, we associate a quadruple
$b_{\sigma ,\pm }=(b_1,b_2,b_3,b_4)$ to the pairs
$(q_i,a_2)$ and $(a_{n-1},q_j^{-1})$  in the reduced code
$c_r(\sigma )$ as follows
$$
b_{(\sigma ,+)}=
\cases
 (o_{i-1},l_i,o_i^{-1},h_i^{-1}) & \qquad \qquad  \text{if}\, \,
a_2=h_{i-1}^{-1}; \\
(l_i,o_i^{-1},h_i^{-1},o_{i-1}) & \qquad \qquad  \text{if}\, \,
a_2=l_{i-1}^{-1}; \\
 (o_i^{-1},h_i^{-1},o_{i-1},l_i) &  \qquad \qquad  \text{if}\, \,
a_2=l_{i+1}; \\
 (h_i^{-1},o_{i-1},l_i,o_i^{-1}) & \qquad \qquad  \text{if}\, \,
a_2=h_{i+1} \\
\endcases
$$
and
$$
b_{(\sigma ,-)}=
\cases
 (o_{j-1},l_j,o_j^{-1},h_j^{-1}) & \qquad \qquad  \text{if}\, \,
a_{n-1}=h_{j-1}; \\
 (l_j,o_j^{-1},h_j^{-1},o_{j-1}) & \qquad \qquad  \text{if}\, \,
a_{n-1}=l_{j-1};  \\
 (o_j^{-1},h_j^{-1},o_{j-1},l_j ) &  \qquad \qquad  \text{if} \, \,
a_{n-1}=l_{j+1}^{-1}; \\
 (h_j^{-1},o_{j-1},l_j,o_j^{-1}) & \qquad \qquad  \text{if}\, \,
a_{n-1}=h_{j+1}^{-1}. \\
\endcases
$$
We define
$$c_1c_2=(a^{\prime }_1,\dots ,a_{m_1}^{\prime },a^{\prime \prime }_1,\dots ,
a_{m_2}^{\prime \prime })$$
for $c_1=(a^{\prime }_1,\dots ,a_{m_1}^{\prime })$ and
$c_2=(a^{\prime \prime }_1,\dots ,a_{m_2}^{\prime \prime })$, and
$$c_1^{-1}=(a_{m_1}^{\prime -1},\dots ,a_{1}^{\prime -1}).$$
Put
$$c_r(\sigma )_{(k+)}=(a_{k+1},\dots , a_{n-1}),$$
$$c_r(\sigma )_{(k-)}=(a_{2},\dots ,a_{k}),$$
$$r_{\sigma }=(a_2,,\dots , a_{n-1})b_{(\sigma ,-)}
(a_2,,\dots , a_{n-1})^{-1}b_{(\sigma ,+)},$$
$$r_{\sigma ,k}=c_r(\sigma)_{(k+)}b_{(\sigma ,-)}c_r(\sigma )_{(k+)}^{-1}
c_r(\sigma )_{(k-)}^{-1}b_{(\sigma ,+)}c_r(\sigma )_{(k-)},$$
$$R_1(\sigma ,k)=[(c_r(\sigma)_{(k+)}b_{(\sigma ,-)}c_r(\sigma
)_{(k+)}^{-1}),\,
(c_r(\sigma )_{(k-)}^{-1}b_{(\sigma ,+)}c_r(\sigma )_{(k-)})]
$$
and
$$R_{\nu}(\sigma ,k)=r_{\sigma ,k}^{\mu -1}
R_1(\sigma ,k)r_{\sigma ,k}^{1-\mu},
$$
where $\nu =2 \mu $, $c_r(\sigma )=(a_1,a_2,\dots , a_{n-1},a_n)$,
$[c_1,c_2]=c_1c_2c_1^{-1}c_2^{-1}$, and each code
$r_{\sigma ,k}$ can be obtained from $r_{\sigma }$ by means of a cyclic
permutation. It is clear that $r_{\sigma ,k}$ and $R_{\nu }(\sigma ,k)$ are
reduced codes.

Without loss of generality, we can choose $\sigma $ and $\gamma $
such that
$c(\sigma )=c_r(\sigma )$ and $\gamma =\delta \circ h_s\circ C_s\circ h_s^{-1}
\circ \delta ^{-1}$,
where $C_s$ is a circle of a small radius with the center at $q_s$, $h_s$ is
a path along the line $v_1=s$, and $\delta $ is the shortest path along
$\partial D$ connecting $u_0$ and
the intersection point of $\partial D$ and the ray $\{ v_1=s, \, \, v_2>0\} $.
Let  $w_1>w_2> \dots > w_m>0$ be the sequence of values of $v_2$
corresponding to
the intersection points of the ray $\{ v_1=s, \, v_2>0\} $ and $\sigma $,
and let  $k_1,k_2,\dots , k_m$ be the indices of
$a_{k_l}$  in $c_r(\sigma )$ and corresponding to these values of $v_2$
(by definition, these $a_{k_l}$ are equal to $h_s^{\varepsilon _{k_l}}, \, \,
\varepsilon _{k_l}= \pm 1$).
Then
$$\align
 c((\gamma )H(\sigma )^{\nu })= \phantom{(\sigma, k_l)^{-sgn(\sigma )
\varepsilon _{k_l}},
o_{s-1},l_s, o_s^{-1},h_s^{-1},(\prod_{l=1}^{m}R_{\nu}(\sigma, k_l)
^{-sgn(\sigma )\varepsilon _{k_l}})^{-1},u_0^{-1}) } & \\
=(u_0,\prod_{l=1}^{m}
R_{\nu}(\sigma, k_l)^{-sgn(\sigma )\varepsilon _{k_l}},
o_{s-1},l_s, o_s^{-1},h_s^{-1},(\prod_{l=1}^{m}R_{\nu}(\sigma, k_l)
^{-sgn(\sigma )\varepsilon _{k_l}})^{-1},u_0^{-1}) &
\endalign
$$
if $q_s$ is not the  starting point of $\sigma $. If $q_s$ is the  starting
point of
$\sigma $, then the code
$c((\gamma )H(\sigma )^{\nu })$ is equal to
$$(u_0,\prod_{l=1}^{m-1}
R_{\nu}(\sigma, k_l)^{-sgn(\sigma )\varepsilon _{k_l}},
\delta _0(\sigma ),(\prod_{l=1}^{m-1}R_{\nu}(\sigma, k_l)
^{-sgn(\sigma )\varepsilon _{k_l}})^{-1},u_0^{-1})
$$
if $a_2= h_{s-1}^{-1}$, where
$$\delta _0(\sigma )=r_{\sigma}^{\mu -1}\delta ^{\prime}_0(\sigma )
r_{\sigma}^{1-\mu},$$
and
$$\align
\delta ^{\prime}_0(\sigma )= \phantom{abcd a_{n-1})b_{(\sigma -)}
(a_2,\dots , a_{n-1})^{-1}b_{(\sigma +)}
(a_2,\dots , a_{n-1})b_{(\sigma -)}^{-1}
(a_2,\dots , a_{n-1})^{-1})} & \\
=(a_2,\dots , a_{n-1})b_{(\sigma -)}
(a_2,\dots , a_{n-1})^{-1}b_{(\sigma +)}
(a_2,\dots , a_{n-1})b_{(\sigma -)}^{-1}
(a_2,\dots , a_{n-1})^{-1}), &
\endalign
$$
and
$c((\gamma )H(\sigma )^{\nu })$ is equal to
$$(u_0,\prod_{l=1}^{m}
R_{\nu}(\sigma, k_l)^{-sgn(\sigma )\varepsilon _{k_l}},
\delta (\sigma ),(\prod_{l=1}^{m}R_{\nu}(\sigma, k_l)
^{-sgn(\sigma )\varepsilon _{k_l}})^{-1},u_0^{-1})
$$
if $a_2\neq h_{s-1}^{-1}$, where
$$
\delta = \cases (o_{s-1},\delta _0(\sigma ), o_{s-1}^{-1}) & \qquad \text{ if}
\, \, a_2=l_{s-1}^{-1},  \\
(o_{s-1},l_s,\delta _0(\sigma ),l_s^{-1}, o_{s-1}^{-1}) & \qquad
\text{if}\, \,  a_2=l_{s+1}, \\
(o_{s-1},l_s,o_s^{-1}\delta _0(\sigma ),o_s,l_s^{-1},
o_{s-1}^{-1}) & \qquad  \text{if}\, \,
a_2=h_{s+1}.
\endcases
$$

Without loss of generality, we can assume that  $\nu \gg 2$. In addition, it
is easy to see that $c_r(r_{\sigma ,k_i}r_{\sigma ,k_j}^{-1})\neq \emptyset $
for $k_i\neq k_j$, since $c_r(\sigma )$ is the reduced code.
Therefore, $c_r((\gamma )H(\sigma )^{\nu })\neq c_r(\gamma )=
(u_0, o_{s-1},l_s, o_s^{-1},h_s^{-1},
u_0^{-1})$. Lemma 5.4 is proved.
\edm
\edm

The diffeomorphism $h$ representing the element $G\in C(X)\subset B_p[D,K]$
is defined up to isotopy,
hence by Lemma 5.4 we can assume that $h(\sigma _1)=
\sigma _1$, $h(q_i)=q_i$ for $i=1,\, 2$.

Let us show that there exists a diffeomorphism
$g_2$ representing an element $G_2\in C_1(X)$ such that
$g_2\circ h$ leaves fixed the paths $\sigma _1$ and $\sigma _2$.
In fact, consider the path $h(\sigma _2)$. The point $q_2$ is the  starting
point of this path and some point $q_r\neq q_1, q_2$ is the end point.
Let $s:[0,1]\to D$ be a parametrization of this path such that $s(0)=q_r$.
Consider the reduced code $c_r(h(\sigma _2))=(a_1,\dots , a_n)$ of the path
$h(\sigma _2)$, $a_1=q_r$, $a_n=q_2^{-1}$.
Let for some $i$ the symbol $a_i$ be equal to $o_j^{\pm 1}$, where $1<j<p$.
We choose among all such  $i$ an index $i_0$ such that
$a_{i_0}=o_{j_0}^{\pm 1}$ for which the following conditions are fulfiled:
\roster
\item"{(i)}" if $j_0=2$, then there is no other intersection point of
$h(\sigma _2)$ and $\sigma _2$ lying in $\sigma _2$ between $s(t_{i_0})$ and
$q_3$\, ;
\item"{(ii)}" if $j_0=r$, then there is no other intersection point of
$h(\sigma _2)$ and $\sigma _r$ lying in $\sigma _r$ between $s(t_{i_0})$ and
$q_r$\, ;
\item"{(iii)}" if $j_0=r+1$, then there is no other intersection point of
$h(\sigma _2)$ Ï $\sigma _{r+1}$ lying in $\sigma _{r+1}$ between $s(t_{i_0})$
and $q_{r+1}$;
\item"{(iv)}" if $j_0\neq 2,r,r+1$, then there is no other intersection
point of $\sigma _{j_0}$ and $h(\sigma _2)$ lying in $\sigma _{j_0}$ between
either $q_{j_0}$ and $s(t_{i_0})$ or $s(t_{i_0})$ and $q_{j_0+1}$.
\endroster
Consider one of these cases (the other cases are similar). For instance, let
$j_0\neq 2,r,r+1$ and assume that there is no
intersection  point of $\sigma _{j_0}$ and $h(\sigma _2)$ lying on
$\sigma _{j_0}$ between  $q_{j_0}$ and $s(t_{i_0})$.
Denote by $\widetilde{h(\sigma _2)}$ a path consisting of a part of
$h(\sigma _2)$ starting at $q_r$ and ending at $s(t_{i_0})$ and a part of
 $\sigma _{j_0}$ starting at $s(t_{i_0})$ and ending at $q_{j_0}$. Let us
choose a smooth path  $\widetilde{\sigma }$, sufficiently closed to
$\widetilde{h(\sigma _2)}$ such that
$\widetilde{\sigma }$ connects $q_r$ and
$q_{j_0}$, $\widetilde{h(\sigma _2)}\cap \widetilde{\sigma }=
\{ q_r, q_{j_0}\}$, and such that if we move along $\widetilde{\sigma }$
starting at $q_r$, then the path $\widetilde{h(\sigma _2)}$ is situated to
the right from $\widetilde{\sigma }$. Perform a half-twist
$H(\widetilde{\sigma })\in C_1(X)$. Let a diffeomorphism $h_1$ be a
representative of $H(\widetilde{\sigma })$. It is easy to see that
$h_1(h(\sigma _2))$
is isotopic to the path having the following code
$c(h_1\circ h)=(\widetilde a_0,\dots ,
\widetilde a_{n-i_0-l})$, where $l$ is a non-negative integer,
$\widetilde a_0=q_{j_0}$, $\widetilde a_j=a_{j+i_0+l}$.

Denoting again by $h$ the diffeomorphism $h_1\circ h$ and repeating
the process described above, we can assume that the code $c(h(\sigma _2))=
(a_0,\dots , a_n)$ of the curve $h(\sigma _2)$ satisfies the condition:
for any $i$ the symbol $a_i\neq  o_j^{\pm 1}$, where $0<j<p$.
We can assume for
definiteness that $a_1=l_{r+1}$ and $a_{n-1}=l_{0}$ (the other cases are
similar).
In this case, it it is easy to see that
$$c(h)=(q_r,l_{r+1}\dots ,l_{p-1},o_p^{-1},h_{p-1}^{-1},\dots , h_1^{-1},o^0,
l_0, q_2^{-1}).$$
It is easy to check that the path $(h(\sigma _2))\Dl _{r,p}\circ
\Dl _{3,p}\circ \Dl _{1,2}^2$ is isotopic to $\sigma _2$. But
$\Dl _{r,p}\circ \Dl _{3,p}\circ \Dl _{1,2}^2\in C_2(X)\subset C_1(X)$. Thus,
multiplying $G$
by $\Dl _{r,p}\circ \Dl _{3,p}\circ \Dl _{1,2}^{2}$, we can assume that the
diffeomorphism $h$ representing $G\in C(X)$ leaves fixed the path $L_{1,3}$.

Repeating the  stated above consecutively for $\sigma _3, \dots ,
\sigma _{p-1}, \sigma _0$, we can assume that the diffeomorphism $h$
representing $G\in C(X)$ leaves fixed the path $L_{0,p}=
\sigma _0\cup L_{1,p}$. In this case considering the code of the path
$h(\sigma _p)$, it is easy to see that $h(\sigma _p)$ is isotopic to
$\sigma _p$ in $(D,K)$. Therefore we can assume that $h$ leaves fixed
the diameter $L_{0,p+1}=\sigma _0\cup L_{1,p}\cup \sigma _p$. But in this
case, $h$ is a representative of the identity element of $B_p[D,K]$.
Thus, $C_1(X)=C(X)$.

To prove that $C_1(X)=C_2(X)$, it is sufficient to show that
if a half-twist $H(\sigma )$ is given by a simple path $\sigma $
starting and ending at $K$ and non-intersecting with $\sigma _1$, then
$H(\sigma )\in C_2(X)$.

We note that, for any $G\in B_p[D,K]$, the element $G^{-1}H(\sigma )G$ is
a half-twist given by the path $(\sigma )G$. Therefore, to prove that
any such half-twist $H(\sigma)\in C_2(X)$, one repeats the arguments
stated above using induction on $l_c(\sigma )$.
\edm

\subheading{\S 6.\ Smooth isotopy of fiber space}

Let $M$ be a smooth variety. By definition, a diffeomorphism $F :M\times
[0,1]\to M\times [0,1]$ (or simply $F_t:M\to M$) is a smooth isotopy if
\roster
\item $F(M\times \{ t \} )= M\times \{ t \}$ for all $t\in [0,1]$\, ;
\item $F_0=F_{|M\times \{ 0 \} } :M\times \{ 0\} \to M\times \{ 0 \}$ is
the identity map.
\endroster
Without loss of generality, we shall assume that the isotopy $F_t$ satisfies
an additional condition
\roster
\item"{(3)}" $F_t=F_0$ if $t\leq \varepsilon $ and $F_t=F_1$  if
$t\geq 1-\varepsilon $ for some $\varepsilon >0$.
\endroster
Indeed, instead of a smooth isotopy $F_t,$ we can consider a smooth isotopy
$\widetilde F_t=F_{h(t)}$, where $h:[0,1]\to [0,1]$ is a smooth monotone
function such that $h(t)=0$ if $t\leq \varepsilon $ and $h(t)=1$ if
$t\geq 1- \varepsilon $.

By definition, the composition of smooth isotopies $F'_t$ and
$F^{\prime \prime}_t$ is the smooth isotopy
$F_t=F^{\prime \prime}_t\circ F^{\prime}_t$ given by $F_t =F'_{2t} $ if
$t\leq \frac{1}{2}$ and
$F_t=F^{\prime \prime}_{2t-1}\circ F'_1$ if $t\geq \frac{1}{2}$\, .

Let $U$ be a neighbourhood in $M$ and $\partial U$ its boundary.
Let $F_t:\overline U\to \overline U$ be a smooth isotopy.  It is evident
that if the restriction of $F_t:\overline U\to \overline U$ to a
neighbourhood of $\partial U$ is the identity map for all $t\in [0,1]$, then
$F_t$ can be extended to a smooth isotopy $\widetilde F_t:M\to M$ such that
$\widetilde F_{t|M\setminus U}$ is the identity map for all $t$.

Let $B_1$ and $B_2$ be two plane curves. As in Section 1,  let

$K_i(x)=\{y\bigm| (x,y)\in B_i\}$\, ,
$i=1,2$\, , $(K_i(x)=$ projection to $y$-axis of $\pi^{-1}(x)\cap B_i)$,

$$N_i=\{x\bigm| \# K_i(x)\lvertneqq p\}.$$
$$M_i'=\{(x,y)\in B_i\bigm| \pi\bigm |_{B_i} \, \text{is not \'etale at} \,
(x,y) \}\
(\pi(M'_i)=N_i).$$

Let $E_R$ (resp. $D_R$) be a closed disk of radius $R$ with center at the
origin $o$ on $x$-axis (resp. $y$-axis) such that
$M_i'\subset E_R\times D_R, $\ $N_i\subset\Int(E).$

Assume that $\# K_i(o)=p$.

For each $u_{i,j}\in N_i$, $j=1,\dots , n$, we choose a disc $E_{i,j}$
of small radius $\epsilon \ll 1$ with center at $u_{i,j}$ and choose
simple paths $T_{i,1},\dots ,T_{i,n}$ connecting $u_{i,j}$ with $o$ such
that $<T_{i,1},\dots ,T_{i,n}>$ is a bush. For each $j=1,\dots ,n$
we choose a small tubular neighbourhood $U_{i.j}$ of $T_{i,j}$ and choose a
disc $E_{o}$ of radius $\epsilon \ll 1$ with center at $o$ such that
\roster
\item $U_{i,j_1}\cap U_{i,j_2}\subset E_o$ for $j_1\ne j_2$\,;
\item $E_{i,j}\cap E_{o}=\emptyset $\, for all $j$\, ;
\item the set
$$U_{\Gamma _i}=(\bigcup\limits^n_{j=1}U_{i,j})\cup (\bigcup\limits^n_{j=1}
E_{i,j})\cup E_o$$
is diffeomorphic to a disc and its boundary $\partial U_{\Gamma _i}$ is a
smooth simple loop.
\endroster
Such $U_{\Gamma _i}$ will be called {\it a tubular neighbourhood} of a
$g$-base $\Gamma _i=(l(T_{i,1}),\dots ,l(T_{i,n}))$ of the fundamental group
$\pi _1(E_R\setminus N_i,o)$.

The following lemmas are well-known.

\proclaim{Lemma 6.1} In the notation described above, for
$N_i\subset E_R$,
$i=1,\, 2$, let
$U_{\Gamma _i}$ be two tubular neighbourhoods of $g$-bases $\Gamma _i$ of
$\pi_1(E_R\setminus K_i,o)$. Assume that $\# N_1= \# N_2$. Then there exists a
smooth isotopy $f_t: E_{R}\to E_R$, $t\in [0,\, 1]$, such that
\roster
\item $f_{t}$ is the identity map in a neighbourhood of the boundary of
$E_R$\, ;
\item $f_t(o)=o$ for all $t\in [0,\, 1]$\,;
\item $f_1(E_{1,j})=E_{2,j}$ for each $j=1,\dots , \, n$\, .
\item $f_1(U_{\Gamma _1})=U_{\Gamma _2}$\, .
\endroster
\endproclaim

\proclaim{Lemma 6.2} The isotopy $f_t$ from Lemma 6.1 can be extended to a
smooth isotopy $F_t:  E_R \times \C ^1 \to E_R  \times \C ^1 $ such that
\roster
\item $F_{t| E_R \times D_R}=f_{t} \times Id$ \, ;
\item $F_t$ is the identity map outside $E_{R} \times D_{2R}$\, .
\endroster
\endproclaim

\proclaim{Lemma 6.3} Let smooth real functions $\alpha (u,v), \beta (u,v)$
satisfy the following inequalities
$$\varepsilon _1+\varepsilon _2<\alpha (u,v)<1-\varepsilon _1-\varepsilon _2;$$
$$\varepsilon _1+\varepsilon _2< \beta (u,v)<
1-\varepsilon _1-\varepsilon _2$$
for all $(u,v)\in E_1=\{\, \sqrt{u^2+v^2} \leq 1 \, \} $ and for some
positive $\varepsilon _1, \varepsilon _2\ll 1$. Then there exists a smooth
real function $f_{(\alpha ,\beta )t}(z,u,v)$, $(z,t,u,v)\in [0,1]\times
[0,1]\times E_1$, satisfying the following conditions:
\roster
\item $f_t(z,u,v)$ is
a monotone function for each fixed $(t,u,v)\in [0,1]\times E_1$ \, ;
\item $f_t(z,u,v)\equiv z$ if $0\leq z \leq \varepsilon _1$;
\item $f_t(z,u,v)\equiv z$ if $1- \varepsilon _1 \leq z \leq 1$;
\item $f_t(z,u,v)=z+t(\beta (u,v)-\alpha (u,v))$  if $\alpha (u,v)-
\varepsilon _2\leq z\leq \alpha (u,v)+\varepsilon _2$.
\endroster
\endproclaim

\proclaim{Proposition 6.4} Let $\Cal B=(b_{1}(x),\dots ,b_{p}(x))$
be a collection of non-intersecting sections of the projection
$\pi : E\times D_R \to E=\{x\in \C^1\, \mid \, |x|\leq 1\} $\, . Then there
exists a smooth isotopy $F_t:E\times D_R \to E\times D_R$ such that
\roster
\item $F_t(x,y)=(x,F_{t,x}(y))$ for all $t$ and $x$\, ;
\item $F_1(\Cal B)$ is a collection of constant sections.
\endroster
Moreover, if all $b_j(x)$ are constant sections (equal to
$b_j$) over a neighbourhood of the boundary of $E$, then the isotopy
$F_t$ can be chosen in such a way that
\roster
\item"{(3)}" $F_t$ is the identity map in the
neighbourhood of the boundary of $E\times D_R$ for $t\in [0,1]$\, ;
\item"{(4)}" $F_1(x,b_j(x))=(x,b_j)$\, .
\endroster
\endproclaim
\demo{Proof} The Proposition follows from the  the following Lemma.

\proclaim{Lemma 6.5} Let $\Cal B=(b_{1}(x),b_2,\dots ,b_{k})$
be a collection of non-intersecting sections of the projection
$\pi : E\times D_R \to E=\{x\in \C^1\, \mid \, |x|\leq 1\} $\, such that
$b_2,\dots , b_k$ are constant sections. Then there exists a smooth isotopy
$F_t:E\times D_R \to E\times D_R$ such that
\roster
\item $F_t(x,y)=(x,F_{t,x}(y))$ for all $t$ and $x$\, ;
\item $F_t(x,b_j)=(x,b_j)$ for all $t$ and $j$\, ;
\item $F_1(x,b_1(x))$ is a constant section.
\endroster
Moreover, if $b_1(x)$ is a constant section (equal to
$b_1$) over a neighbourhood of the boundary of $E$, then the isotopy
$F_t$ can be chosen in such a way that
\roster
\item"{(4)}" $F_t$ is the identity map in a
neighbourhood of the boundary of $E\times D_R$ for $t\in [0,1]$\, ;
\item"{(5)}" $F_1(x,b_j(x))=(x,b_j)$\, .
\endroster
\endproclaim
\demo{Proof}
Consider the set $S=\{ y\in D_R\, \mid \, y=b_1(x),\, x\in E\, \} $\, .
At first, assume that we can find a simply connected neighbourhood
$U\subset D_R$ of $S$ such that $b_j\not\in U$ for $j=2,\dots , \, k$. By
Riemann's Theorem, there exists a complex-analytic bijective morphism
$\varphi :U\to U_1$, where $U_1=\{ \, z=z_1+iz_2\in \C \, \mid \, 0<z_1<1,\,
\, \, 0<z_2<1,\, \} $\, . Since $S$ is a compactum, then there exist
$\varepsilon _1, \varepsilon _2>0$ such that
$$ \varepsilon _1+\varepsilon _2<z_l<1-\varepsilon _1-\varepsilon _2, \,
\, \, \, \, \, \, \, \, l=1,2, $$
for $z_1+iz_2\in \varphi (S)$.
We have two smooth functions $\alpha _1(x),\, \, \alpha _2(x)$, where
$\varphi (b_1(x))=\alpha _1(x)+i\alpha _2(x)$, satisfying
$$\varepsilon _1+\varepsilon _2  < \alpha _l(x)<
1-\varepsilon _1-\varepsilon _2 .$$
By Lemma 6.3, for any two smooth functions $\beta _1(x)$ and $\beta _2(x)$
(for instance, if $\beta _1(x)$ and $\beta _2(x)$ are constant functions
such that $\varphi (b_1)=\beta _1+i\beta _2$) satisfying
$$ \varepsilon _1+\varepsilon _2<\beta_l(x)<1-\varepsilon _1-\varepsilon _2, \,
\, \, \, \, \, \, \, \, l=1,2, $$
there exists a smooth isotopy $F_t :E\times U_1\to E\times U_1$ given by
$$\align
F_t^*(x) & \equiv x\, ; \\
F_t^*(z_1) & =f_{(\alpha _1,\beta _1)th(|x|)}(z_1,u,v)\, ; \\
F_t^*(z_2) & =f_{(\alpha _2,\beta _2)th(|x|)}(z_2,u,v)\, ,
\endalign
$$
where $u+iv=x$, $f_{(\alpha _l,\beta _l)t}(z_l,u,v)$ is a function
described in Lemma 6.3, and $h(r)$ is a smooth monotone real function such that
$h(r)=1$ if $r<\varepsilon $ and $h(r)=0$ if $r>1-\varepsilon $. It is easy
to check that the constructed isotopy $F_t$ satisfies
properties (1) - (4) of Lemma 6.5.

If $S$ is not ``wound" around every point $b_j$, $j=2,\dots , k$, then
we can find such a simply connected neighbourhood $U$ of $S$.
Otherwise we will show that there exists a sequence of isotopies which
``unwinds" $S$. To show this, we fix the point $x_0=1\in E$. Let $\gamma
(t)=x_t$
be a smooth simple (without self-intersections) path joining $x_0$ with a
point $x_1\in E$. Denote by $n_{j,\gamma }(x_1)$ the integral part of the
number of rotations around the point $b_j$ of a point moving along the path
$b_1(\gamma (t))$. Evidently, $n_j(x_1)=n_{j,\gamma }(x_1)$ depends only on
the point $x_1$, and  not   on $\gamma $ joining $x_0$ and $x_1$,
since $E$ is simply connected and $b_1(x)$ is a smooth function. Denote by
$$E_{n_2,\dots ,n_k}=\{ \, x\in E\, \mid \, n_j(x)=n_j,\, \, j=2,\dots ,k\,
 \} $$
and let $E_{\max }$ (resp. $E_{\min }$) be the set $E_{n^0_2,\dots ,n^0_k}$,
where $n^0_j$ is a local maximum (resp. minimum) of $n_j(x)$ for all $j$.
Evidently, if $E_{\max }=E_{\min }=E_{0,\dots ,0}$, then we can find a simply
connected neibourhood $U$ of $S$ such that $b_j\not\in U$ for $j=2,\dots ,k$.

Denote by $\overline l_j$ the line passing through $b_1$ and $b_j$ and choose
a coordinate $u_j$ in $\overline l_j$ such that $u_j(b_1)=0$\, ,
where $b_1=b_1(x_0)$. Let $l_j= A_j y_1+B_jy_2$ be a linear function such that
$\overline l_j=\{ \, l_j(y_1,y_2)=0\, \}$, where $ y_1=\text{Re}\, y$ and
$y_2=\text{Im}\, y$.

Let
$$S_j=S\cap \{ \, y\in \overline l_j\, \mid \, u_j(y)< u_j(b_j) \, \,
\text{if}\, \, u_j(b_j)>0\, \, \text{and}\, \, u_j(y)> u_j(b_j)\, \,
\text{if}\, \, u_j(b_j)<0\,  \} .$$
Denote by $C_j=\{ \, x\in E\, \mid \, y=b_1(x)\in S_j\, \} $\, .
Consider $G=E\setminus \bigcup\limits^k_{j=2}C_j$
which is a disjoint union of a finite number of connected components.
Any two neighbouring $G_m$ and $G_l$ are separated by the connected component
$C_j^0$ of $C_j$ for some $j$. If $n_j(m)=n_j(l)$, where $n_j(m)=n_j(x)$ for
$x\in G_m$, then we change $G_m$ and $G_l$ to the union $G_m\cup G_l\cup
C_j^0$.

Evidently, there exists a connected component $G^0$ of $G$ whose boundary
$\partial G^0$ is connected and therefore,
$G^0$ is simply connected. Moreover, it is easy to see that there exists a
simply connected neighbourhood $U_0$ of $\{ \, y\in D_R\, \mid \, y=b_1(x),
x\in \overline G^0\, \} $ such that $b_j\not \in U_0$ for $j=2,\dots , k$.
There is $j_0$ such that $\partial G^0$ is a subset of $\{ \, l_{j_0}(b_1(x))
=0 \, \}$. By the choice of $G^0$, we can assume that $l_{j_0}(b_1(x))\geq 0$
for all $x\in G^0$. Since $l_{j_0}(b_1(x))$ is a smooth function, then
the set $\{ \, l_{j_0}(b_1(x))=\delta \, \}$ is a smooth curve for $\delta $
close to $0$. Therefore, without loss of generality, we can change $G^0$
by a set $\widetilde G^0$ such that
\roster
\item $G^0$ is closed to $\widetilde G^0$ and contained in $\widetilde G^0$\, ;
\item the boundary of $\widetilde G^0$ is the subset of
$\{ \, l_{j_0}(b_1(x))=\delta \, \}$ for some $\delta <0$ and close to $0$\, ;
\item the closure of $b_1(\widetilde G^0)$ is contained in $U_0$, where $U_0$
is a simply connected open set such that $b_j\not\in U_0$ for $j=2,\dots
,k$\, .
\endroster

As above, we can find holomorphic coordinate $z=z_1+iz_2$, $z=\varphi(y)$, in
$U_0$ such that $U_0\simeq \{\, 0<z_1<1\, \} \times \{\, 0<z_1<1\, \} $.
Evidently, for the closure $G^1$ of $\widetilde G^0$ there exist smooth
functions $\alpha : G^1\to U_0$, $z=\varphi (b_1(x))=\alpha _1(x)+
i\alpha _2(x)$, and $\beta (x)=\beta _1(x)+i\beta _2(x)$
such that $\beta (x)=\varphi (b_1(x))$ for $x$ lying in a neighbourhood of
$\partial G^1$ and $l_{j_0}(\beta (x))<0$ for all $x\in G^1$.
 As above applying Lemma 6.3, we can find a smooth isotopy
$F_t:G^1\times U_0\to G^1\times U_0$ such that $F_1(x,b_1(x))$ does not meet
the set $\{ (x,y)\, \mid \, x\in G^1,\, y\in \cup \overline l_j\}$ and such
that $F_t$ is the identity map in a neighbourhood of the boundary of $G^1
\times U_0$. Hence this isotopy can be extended to the isotopy of $E\times
D_R$. Evidently, for the image $F_1(b_1(x))$ of the section $b_1(x)$, we can
repeat the construction of the set $G$ and observe that the number of connected
components of $G$ is decreased. Hence,  after several similar steps, we
construct the  desired isotopy as the composition of isotopies constructed
in each step.

 From the consideration described above, the following remark follows.
\proclaim{Remark 6.6} \ Let $C\subset D$ be a connected set such that
the closure $\overline C$ and the boundary $\partial D$ of $D$ have a
non-empty intersection. Assume that $b_1(x)\equiv b_1$ for $x\in C$.
Then on each step (except the last one), $G^0$ can be chosen in such a way that
$C\not\subset G^0$.
\endproclaim

The proof of Lemma 6.5 in the case when $b_1(x)$ is a constant section over a
neighbourhood of the boundary of $D$ follows from Remark 6.6.\hfill $\qed$
for Lemma 6.5
\enddemo
\enddemo

\proclaim{Remark 6.7} Let $\Cal B_1=(b_{1,1}(x),\dots ,b_{1,p}(x))$ and
$\Cal B_2=(b_{2,1}(x),\dots ,b_{2,p}(x))$ be two collections of
non-intersecting sections of the projection
$\pi : E\times D_R \to E=\{x\in \C^1\, \mid \, |x|\leq 1\} $ such that
for all $j$ the sections $b_{1,j}(x)$ and $b_{2,j}(x)$ coincide over
a neighbourhood $U$ of the boundary of $E$. It is easy to see that the smooth
isotopies $F'_t:E\times D_R \to E\times D_R$ for $\Cal B_1$ and
$F^{\prime \prime }_t:E\times D_R \to E\times D_R$ for $\Cal B_2$
from Proposition 6.4 can be chosen in such a way that
$F'_t$ and $F^{\prime \prime }_t$ coincide over $U$.
\endproclaim

Let $\Cal B=(b_{1}(x),\dots ,b_{p}(x))$ be a collection of
non-intersecting sections of the projection
$\pi : E\times D_R \to E=\{x\in \C^1\, \mid \, |x|\leq 1\} $. Denote by
$$U_\varepsilon (b_j(x))=\{ (x,y)\in D\times E_R\, \mid \,
|y-b_j(x)|<\varepsilon \, \} $$
a tubular neighbourhood of the section
$b_j(x)$ and put $U_{\varepsilon }(\Cal B)=\cup U_\varepsilon (b_j(x))$.
Let $F_t:E\times D_R \to E\times D_R$ be a smooth isotopy. Consider
a smooth map $\widetilde F'_t:U_{\varepsilon }(\Cal B)\times [0,1]\to
E\times D_R\times [0,1]$ given by
$$\widetilde F'_t(x,y)=(x,y+(F_t(b_j(x))-b_j(x))$$
if $(x,y)\in U_\varepsilon (b_j(x))$. We observe that $\widetilde
F'_{t|\Cal B}=
F_{t|\Cal B}$.

Using the standard technique of pasting together vector fields, one
can prove the following
\proclaim{Lemma 6.8} Let an isotopy $F_t:E\times D_R \to E\times D_R$
have properties $(1) - (4)$ of Proposition 6.4. Then for some
$\varepsilon _1\ll \varepsilon $ the map
$\widetilde F'_t:U_{\varepsilon }(\Cal B)\times [0,1]\to E\times D_R\times
[0,1]$ can be extended to an isotopy $\widetilde F_t:E\times D_R \to E
\times D_R $ having properties $(1) - (4)$ of Proposition 6.4.
\endproclaim

\proclaim{Lemma 6.9} Let $\Cal B_i=(b_{i,1}(x),\dots ,b_{i,p}(x))$
be two collections of non-intersecting sections of the projection
$\pi : E\times D_R \to E=\{x\in \C^1\, \mid \, |x|\leq 1\} $\, . Then there
exists a smooth isotopy $F_t:E\times D_R \to E\times D_R$ such that
\roster
\item $F_t(x,y)=(x,F_{t,x}(y))$ for all $t$ and $x$\, ;
\item for each $j$ there exists a neighbourhood
$U_j$ of the section $b_j(x)$ such that $F_{1|U_j}$ is
holomorphic in $y$\, ;
\item $F_1(\Cal B_1)=\Cal B_2$ over $E_{R_1}\subset E$ for some $R_1<1$.
\item $F_t$ is the identity map in the
neighbourhood of the boundary of $E\times D_R$ for $t\in [0,1]$\, .
\endroster
\endproclaim
\demo{Proof} By Proposition 6.4 and Lemma 6.8, there exists a smooth isotopy
$\widetilde F_t :E\times D_R \to E\times D_R$ having properties (1) - (2)
of Lemma 6.9 and such that $\widetilde F_1(\Cal B_1)=\Cal B_2$. Let $h(r)$ be
a smooth monotone function such that $h(r)=1$ if $r\leq R_1<1-\varepsilon$
and  $h(r)=0$ if $r\geq 1-\varepsilon$. Then $F_t(x,y)=\widetilde F_{h(|x|)t}
(x,y)$ has all properties of Lemma 6.9.
\enddemo

\proclaim{Lemma 6.10} Let $\Cal B_i=(b_{i,1}(x),\dots ,b_{i,p}(x))$
be two collections of non-intersecting sections of the projection
$\pi : E\times D_R \to E=\{x\in \C^1\, \mid \, x_1=\text{Re}\, x\in [0,1],\, \,
x_2=\text{Im}\, x\in [0,1] \, \} $ such that over $E_{\varepsilon }=
\{ 0\leq x_1\leq \varepsilon \}\times \{ 0\leq x_2\leq 1 \}\cap
\{ 1-\varepsilon \leq x_1\leq 1 \}\times \{ 0\leq x_2\leq 1 \}$ the
sections $(b_{1,1}(x),\dots ,b_{1,p}(x))=(b_{2,1}(x),\dots ,b_{2,p}(x))=
(b_{1},\dots ,b_{p})$ are coinciding constant sections. Let
the geometric braids
$$\overline{ \Cal B}_1=(b_{1,1}(x_1),\dots ,b_{1,p}(x_1))\, \, \, \,
\text{and}\, \, \, \,
\overline{ \Cal B}_2=(b_{2,1}(x_1),\dots ,b_{2,p}(x_1)), \, \, \, \, \,
x_1\in[0,1],$$
are two
representatives of the same element in the braid group $B_p$. Then there
exists a smooth isotopy $F_t:E\times D_R \to E\times D_R$ such that
\roster
\item $F_t(x,y)=(x,F_{t,x}(y))$ for all $t$ and $x$\, ;
\item for each $j$ there exists a neighbourhood $U_j$ of the section
$b_j(x)$ such that $F_{1|U_j}$ is holomorphic in $y$\, ;
\item $F_1(\Cal B_1)=\Cal B_2$ over $\{0\leq x_1\leq 1\}\times \{ \,
|x_2|<\varepsilon _2\, \} $ for some $\varepsilon _2 >0$.
\item $F_1$ is the identity map in the neighbourhood of the boundary of
$E\times D_R$.
\endroster
\endproclaim
\demo{Proof} By Proposition 6.4, there exists a smooth isotopy
$\widetilde F_t$ such that $\widetilde F_1(b_{1,j}(x))=b_{j}$,
$j=1,\dots ,p$, are constant sections.

Fix the point $x_0=(0,0)$. As above, for each section
$\widetilde F_1(b_{2,j_0}(x))$ we can define a function $n_{j_0,j}(x)$
equals to the number of rotations around $b_j$. Evidently,
$\widetilde F_1(\overline{\Cal B}_1)$ and
$\widetilde F_1(\overline{\Cal B}_2)$ are also the
representatives of the same element in $B_p$. Therefore, for each section
$\widetilde F_1(b_{2,j_0}(x))$ the number of rotations $n_{j_0,j}(x)=0$ for
$x\in \{ 1-\varepsilon \leq x_1\leq 1 \}\times \{ 0\leq x_2\leq 1 \}$.
Therefore, Lemma 6.10 follows from Remark 6.6 and Lemma 6.9.
\enddemo

\subheading{\S 7.\ Braid monodromy factorization types and
diffeomorphisms of pairs}

Consider a linear projection $\pi :\CP ^2\to \CP ^1$ with center at
$z\in \CP^2$.

\definition{Definition}\ \underbar{Semi-algebraic curve (with respect to
$\pi $)}

A closed subset $B\subset \CP ^2$, $z\not\in B$, is called a semi-algebraic
curve with respect to
$\pi ,$ if for each point $x\in B$ there exist a neighbourhood $U\subset
\CP ^2$
of $x$ and local analytic coordinates $(z_1,z_2)$ in $U$ such that
\roster
\item $\pi _{|U}$ is given by $\pi (z_1,z_2)= z_1$;
\item  either $B\cap U$ is a smooth section of $\pi _{|U}$ over $\pi (U)$,
or $B\cap U$  coincides with a set given by the  equation $f(z_1,z_2)=0$,
where $f(z_1,z_2)$ is an analytic function.
\endroster

A semi-algebraic curve is (generalized) {\it cuspidal} if in $(2)$  the
function $f$ coincides with $f=z_1^k-z_2^2$, $k\in \Bbb N$. The point
$(0,0)$ in a neighbourhood of which $B$ is given by $z_1^k-z_2^2=0$ is called
a {\it singular} point of $B$.
 \enddefinition

Obviously, any algebraic curve $B\subset \CP ^2$ is a semi-algebraic curve
with respect to a generic projection.

As in the algebraic case, we can resolve the singular points of a
semi-algebraic
curve $B$ by means of a composition of $\sigma$-processes
$\nu :\overline{\Bbb P}^2\to \Bbb P^2$ (we needn't resolve the tangent points
of $B$ with the fibers of
$\pi $) and obtain a non-singular
Rimannian surface $\overline B\subset \overline{\Bbb P}^2$. The
composition $(\pi \circ \nu )_{|\overline B}
:\overline B\to \Bbb P^1$ allows us to introduce a complex structure on
$\overline B$ such that $(\pi \circ \nu )_{|\overline B}$ is a holomorphic map,
but $\nu _{|\overline B} :\overline B\to \Bbb P^2$ is not holomorphic and
it is a $C^{\infty }$-map only.

It is clear that as in the algebraic case one can define the braid monodromy
with respect to $\pi $, braid monodromy factorization, and braid monodromy
factorization type for any semi-algebraic curve.

Let $\sigma : \F  \to \CP ^2$ be $\sigma $-process with center at $z$,
$L=\sigma ^{-1}(z)$.
Denote again by  $\pi :\F  \to \CP ^1$ the composition $\pi \circ \sigma $.

\proclaim{Theorem 7.1} \ Let two (generalized) cuspidal semi-algebraic
curves $B_1$ and $B_2$
have the same braid monodromy factorization type $\Delta (B_1)=\Delta
(B_2)$. Then
there exists a smooth isotopy $F_t:\F \to \F$ such that
\roster
\item  $F_{t|U}$ is the identity map for all $t\in [0,1]$, where $U$ is a
neighbourhood of the exceptional curve $L$;
\item  for each point $p\in B_1$ there exist neighbourhoods $U_1$ of $p$ and
$U_2=F_1(U_1)$ of $F_1(p)$ with local complex coordinates $(x_i,y_i)$ in $U_i$
such that $\pi _{|U_i}$ is given by $(x_i,y_i)\mapsto x_i$ and
$F_t^*(y_2)=\phi (x_1,y_1)$ is a smooth complex function holomorphic in
$y_1$\, ;
\item  for each singular point $s\in B_1$ there exists a neighbourhood
$U\subset \F$ of $s$ such that $F_{1|U}:U\to F_1(U)$ is a holomorphic map;
\item $F_1(B_1)= (B_2)$.
\endroster
\endproclaim

\proclaim{Corollary 7.2} \ Let two (generalized) cuspidal semi-algebraic curves
$B_1$ and $B_2$ have the same braid monodromy factorization type $\Delta (B_1)
=\Delta (B_2)$. Then there exists a diffeomorphism of pairs $F: (\CP ^2, B_1)
\to (\CP ^2, B_2)$ having properties $(2)$ and $(3)$ of Theorem $
7.1.$\endproclaim

\demo{Proof} \  An isotopy $F_t:\F\to \F$ is called {\it compatible} with
a semi-algebraic curve $B$ (with respect to $\pi $) if $F_t(B)$
is a semi-algebraic curve $B$ with respect to $\pi $ for each $t\in [0,1]$.

The required isotopy $F_t$ will be obtained as a composition of a
sequence of smooth isotopies compatible with $B_1$.

Without loss of generality, we can assume that $B_1$ and
$B_2$ are embedded into the same $\F $ and they are (generalized) cuspidal
semi-algebraic curves of degree $p$ with respect to $\pi $. We fix a point
$\infty \in \CP ^1$ such that $\pi ^{-1}(\infty )$ is a generic fiber of
$\pi $ with respect to each $B_1$ and $B_2$. Denote by $\C ^1 = \CP ^1
\setminus \{ \infty \} $ and $\C ^2 = \F \setminus (\pi ^{-1}(\C ^1 )
\cup L)$. We choose coordinates $(x,y)$ in $\C ^2$ such that
$\pi: \C^2\ri \C$ is a projection on the first coordinate.

As in Section 6, let $K_i(x)=\{y\bigm| (x,y)\in B_i\}$\, ,
$i=1,2$\, , $(K_i(x)=$ projection to $y$-axis of $\pi^{-1}(x)\cap B_i)$,
$N_i=\{x\bigm| \# K_i(x)\lvertneqq p\}$,
$$M_i'=\{ (x,y)\in B_i\bigm| \pi\bigm |_{B_i}\,  \text{is not \'etale at} \,
(x,y) \}\
(\pi(M'_i)=N_i).$$

Let $E_R$ (resp. $D_R$) be a closed disk of radius $R$ with center at the
origin $o$ on $x$-axis (resp. $y$-axis) such that
$M_i'\subset E_R\times D_R, $\ $N_i\subset\Int(E).$

Assume that $\# K_i(o)=p$.

For each $u_{i,j}\in N_i$, $j=1,\dots , n$, we choose a disc $E_{i,j}$
of small radius $\epsilon \ll 1$ with center at $u_{i,j}$ and choose
simple paths $T_{i,1},\dots ,T_{i,n}$ connecting $u_{i,j}$ with $o$ such
that $<T_{i,1},\dots ,T_{i,n}>$ is a bush. For each $j=1,\dots ,n$
we choose a small tubular neighbourhood $U_{i.j}$ of $T_{i,j}$ and choose a
disc $E_{o}$ of radius $\epsilon \ll 1$ with center at $o$ such that
$$U_{\Gamma _i}=(\bigcup\limits^n_{j=1}U_{i,j})\cup (\bigcup\limits^n_{j=1}
E_{i,j})\cup E_o$$
is  a tubular neighbourhood of a
$g$-base $\Gamma _i=(l(T_{i,1}),\dots ,l(T_{i,n}))$ of the fundamental group
$\pi _1(E_R\setminus N_i,o)$.

By Lemmas 6.1 and 6.2, we can assume that $E_{i,j}=\{\, x\in \Bbb C^1\,
|\, |x-x_{i,j}|<2\}$, where $x_{i,j}$ is the coordinate of $u_{i,j}$.

Without loss of generality, we can assume that $T_{i,j}\cap E_{i,j}$ is a
radius in $E_{i,j}$. We extend this radius to the diameter $d_{i,j}$ and
denote by $\widetilde T_{i,j}$ the path continuing $T_{i,j}$ along this
diameter, $T_{i,j}\subset \widetilde T_{i,j}$.

{\it Step I.} Since $B_1$ and $B_2$ have the same braid monodromy
factorization type, then by Lemma 3.1 we can choose (and fix) $g$-bases
$\Gamma _i$, $i=1,2$\, , such that the braid monodromy
factorizations $\Dl (B_1)$ and $\Dl (B_2)$ associated to them are
conjugation equivalent.

\proclaim{Lemma 7.3} Let $B\subset \F$ be a semi-algebraic curve and let
$F_t:\F \to \F $ be a smooth isotopy compatible with $B$ such that
$F_t$ is the identity map over the complement of a neighbourhood
$U\subset \CP ^1$, $o\not\in U$.
Then $B$ and $F_1(B)$ have the same braid monodromy factorization.
\endproclaim
\demo{Proof} A family of homomorphisms
$\vp _t:\pi_1(E\setminus N,o)\ri B_p[\C_o, K]$ induced by $F_t$ is
continuous. Therefore, the braid
monodromy
factorizations corresponding to $F_t(B)$ do not depend on $t$, since
$B_p[\C_o, K]$ is a discrete group.
\enddemo

By Lemmas 6.1, 6.2, and 7.3, there exists a smooth isotopy $F_t :\F \to \F$
having properties $(1)-(3)$ of Theorem 7.1 and such that
\roster
\item $F_1(U_{\Gamma _1 }\times \pi ^{-1}(U_{\Gamma _1 }))=U_{\Gamma _2}
\times \pi ^{-1}(U_{\Gamma _2 })$\, ;
\item $F_1(N_1\times \pi ^{-1}(N_1))=N_2\times \pi ^{-1}(N_2)$\, ;
\item the types of singular points $s_{1,j}\in F_1(B_1)$ and $s_{2,j}\in B_2$
coincide over the point $x_{2,j}=x_{j}\in N_2=N$;
\item $B_1$ and $F_1(B_1)$ have the same (up to conjugation equivalence)
braid monodromy factorization.
\endroster

Denote again by $B_1$ its image $F_1(B_1)$.

{\it Step II. } By Lemma 6.9, there exist smooth isotopies $F'_t:\F \to \F$ and
$F^{\prime \prime }_t:\F \to \F$ having properties $(1)$ - $(3)$ of
Theorem 7.1, and such that $F'_1(B_1)$ and $F^{\prime \prime }_1(B_2)$ coincide
over $E_o$, where $E_o$ is a disc with center at the origin in $\C ^1$.
Moreover, we can assume that $F'_1(B_1)$ and $F^{\prime \prime }_1(B_2)$ are
constant sections over $E_o$:
$$K(u)=\{ y_j\in D_R  \, \mid \, y_j=2j-1,\, j\in \Bbb Z,\, 0\leq j
\leq p-1\, \} $$
for all $u\in E_o$.

Denote again by $B_1$ and $B_2$ their images $F'_1(B_1)$ and
$F^{\prime \prime }_1(B_2)$, respectively.
\proclaim{Remark 7.4} The isotopies $F'_t$ and $F^{\prime \prime}_t$ change
braid monodromy factorizations associated to $B_1$ and $B_2$ to conjugation
equivalent factorizations, but they
do not change their braid monodromy factorization types.
\endproclaim

We fix a frame $(\sigma _1,\dots , \sigma _{p-1})$ of $B_p[\Bbb C^1_o,K(o)]$,
$$\sigma _j =[2j-3,2j-1]=\{ \, y\in \Bbb C^1\, \, | \, \,
2j-3\leq \text{Re}\, y\leq 2j-1,\, \text{Im}\, y=0\, \} .$$

{\it Step III. } In $K_i(u_j)=\{ q_{i,1},\dots , \, q_{i,p-1} \}$ let the point
$q_{i,1}$ be the singular point of $B_i$, where $u_j$ is the center of the
disc $E_{j}=E_{1,j}=E_{2,j}$ defined in the definition of tubular
neighbourhood of $g$-base.

\proclaim{Lemma 7.5} Let $B=\{ f(x,y)=0 \} $ be a germ of an analytic curve in
$U=E_{\varepsilon}\times D_{\varepsilon}$ and let $(0,0)$ be a singular point
of $B$ of multiplicity 2 in direction $x=const$ (that is, $\# (B\cap
(\{ u\}\times D_1))=2$ for each $u\in E_{\varepsilon}$). Then there exists a
smooth isotopy $F_t: U\to U$ such that
\roster
\item for all $t\in [0,1]$, $F_t=Id$ in a neighbourhood of the boundary of
$U$\, ;
\item there exists $\varepsilon _1 \ll \varepsilon $ such that $F_{1}(B)\cap V=
\{ y^2 -x^k=0\} $,  where $V=E_{\varepsilon _1}\times D_{\varepsilon _1}$ ;
\item $F_{t|V}$ is holomorphic for each $t$\, .
\endroster
\endproclaim
\demo{Proof} \  By Weierstrass Preparation Theorem we can assume that $B$
is given in $U$ by
$$y^2+h_1(x)y+h_2(x)=0, \tag 1 $$
where $h_i(x)$ are analytic functions. Write $(1)$ in the form
$$(y+\frac{1}{2}h_1(x))^2-(\frac{1}{4}h_1^2(x)-h_2(x))=(y+g_1(x))^2-x^kg_2(x))=0
,$$
where $g_2(0)=c=re^{i\varphi }\neq 0$.

Let $F'_t$ be an isotopy given by $F'_t(x,y)=(x,y+th(|x|)g_1(x))$, where
$h(r)$ is a smooth monotone function such that $h(r)=1$ if $r<\varepsilon _1
\ll \varepsilon $ and $h(r)=0$ if $r>\varepsilon -\varepsilon _1$,

One can show that a smooth map $\widetilde F^{\prime \prime }_t:
E_{\varepsilon _1}\times [0,1]\to E_{\varepsilon }\times [0,1]$ given by
$$\widetilde F^{\prime \prime }_t(x)=x((1+(r-1)t)e^{it\varphi } +t
(g_2(x)-c))^{1/k} $$
can be extended to a smooth isotopy $\widetilde F^{\prime \prime }_t:
E_{\varepsilon }\times[0,1]\to E_{\varepsilon }\times [0,1]$ such that
$\widetilde F^{\prime \prime }_t$ is the identity map in a neighbourhood of
the boundary of $E_{\varepsilon }$. Then the composition $F_t=F^{\prime
\prime }_t
\circ F^{\prime }_t$, where $F^{\prime \prime }_t$ is given in $U$ by
$$ F^{\prime \prime }_t(x,y)=(\widetilde F^{\prime \prime }_t(x),y),$$
satisfies the conditions of Lemma 7.5.
\enddemo

By Lemmas 6.9 and 7.5, there exist smooth isotopies $F^{\prime }_t:\F \to \F$
and $F^{\prime \prime }_t:\F \to \F$
having properties $(1)$ - $(3)$ of Theorem 7.1 and such that
\roster
\item $F'_1(B_1)$ Æ $F^{\prime \prime}_1(B_2)$ coincide over
$\bigcup\limits^n_{j=1}E'_{j}$, where
$E'_{j}\subset E_{1,j}=E_{2,j}=E_{j}$ are some small neighbourhoods of
$u_{i,j}=u_j$;
\item by Lemmas 6.1 and 6.2, we can assume that $E'_{j}=\{\, x\in \Bbb C^1\, \,
|\, \, |x-x_{j}|<2\, \}$, where $x_{j}$ is the coordinate of the point $u_{j}$;
\item in a neighbourhood of $(x_j,0)$ the curves $F'_1(B_1)$ and
$F^{\prime \prime}_1(B_2)$ are given by the equation $y^2=(x-x_{j})^k$\, ;
\item all other $p-2$ branches of $F'_1(B_1)$ (resp.
$F^{\prime \prime}_1(B_2)$) are constant sections over $E'_j$ and
$$K(u'_j)=\{ \, y_1=-1,\dots , y_j=2j-3,\dots ,y_p=2p-3\, \} .$$
for $u'_j=\{ x'_j=x_j+1\}$.\endroster
 Denote again by $B_1$ and $B_2$ their images $F'_1(B_1)$ and
$F^{\prime \prime}_1(B_2)$ respectively.
Without loss of generality, we can assume that $u'_j\in T_j$. Denote by
$u^{\prime \prime}_j$ a point lying in the diameter $\widetilde T_j\cap E_j$
such that $u^{\prime \prime}_j$ is symmetric to $u'_j$ with respect to the
center $u_j$. Let $d'_j$
be a part of the diameter connecting $u'_j$ and $u^{\prime \prime}_j$.

{\it Step IV. } By Lemma 7.3 and Remark 7.4 the isotopies described above
do not
change the braid monodromy factorization types of $B_1$ and $B_2$.
Write the braid monodromy factorizations of the curves $B_1$ and $B_2$
associated to the $g$-base $\Gamma $ fixed above:
$$\align
\Delta ^2  & = \prod\limits_{j=1}^{n} Q_j^{-1}\sigma _1^{\nu _j}Q_j \, \,
\, \,
\text{for}\, \, B_1 ; \\
\Delta ^2  & = \prod\limits_{j=1}^{n} Q^{-1}Q_j^{-1}\sigma _1^{\nu _j}Q_jQ
\, \, \, \,
\text{for}\, \, B_2.
\endalign
$$

We show that in our case there exists a smooth isotopy $F_t:\F \to \F$
having properties $(1)$ - $(3)$ of Theorem 7.1 and such that
\roster
\item $F_t$ is the identity map over the complement of the neighbourhood
$U_{o}$ of the point $o$\, ;
\item $F_1(B_1)=B_2$ over a neighbourhood $U'_{o}\subset U_o$\, ;
\item $F_1(B_1)$ and $B_2$ have the same braid monodromy factorization.
\endroster

In fact, let
$$K=\{ y_j\in \C^1\, \mid \, y_j=2j-3,\, j=1,2,\dots , p\, \} .$$
The half-twist $H_j=H(\sigma _j )\in B_p=B_p[\Bbb C^1_o,K(o)]$ can be
represented by a geometric braid $\overline \sigma _j(s_1)$ in
$\C ^1\times [0,1]$\, :
$$\align \overline \sigma _{j,l}(s_1) & =l\, \, \,  \text{for} \, \,
l=1,\dots, j-1, \, j+2,\dots , p\, ; \\
\overline \sigma _{j,j}(s_1) & =e^{\pi (\beta (s_1)+1)i}+2j-2\, ; \\
\overline \sigma _{j,j+1}(s_1) & =e^{\pi \beta (s_1)i}+2j-2\, , \endalign
$$
where $s_1\in [0,1]$, $\beta(s_1)$\, is real smooth monotone
function such that $\beta (s_1)=0$ for $s_1\in [0,\frac{1}{3}]$ and
$\beta (s_1)=1$ for $s_1\geq \frac{2}{3}$.

The element $H_j^{-1}$ can be represented by a geometric braid
$\overline \sigma _j^{-1}(s_1)$ in $\C ^1\times [0,1]$\, :
$$\align \overline \sigma _{j,l}^{-1}(s_1) & =l\, \, \,  \text{for} \, \,
l=1,\dots, j-1, \, j+2,\dots , p\, ; \\
\overline \sigma _{j,j}^{-1}(s_1) & =e^{\pi (-\beta (s_1)+1)i}+2j-2\, ; \\
\overline \sigma _{j,j+1}^{-1}(s_1) & =e^{-\pi \beta (s_1)i}+2j-2\, , \endalign
$$
The product $Q =H_{j_1}^{\delta _1}\cdot \dots \cdot H_{j_k}^{\delta _k}$,
where $\delta _l=\pm 1$, can be represented by the geometric braid
$\overline Q (s_1)$ in $\C ^1\times [0,k]$\, :
$$\overline Q_l(s_1)= \overline \sigma _{j_m,l}^{\delta _{j_m}}(s_1-m+1)\,
\, \, \text{for}\, \, s_1\in [m-1,m]. $$

Let $U'\subset U_o$ be a neighbourhood of $o$ for which there exists a
diffeomorphism $\varphi : U'\to V=(-1,2k+1)\times (0,2)$, $\varphi (o)=(0,0)$.
Obviously, the paths $T_j$ representing the bush can be chosen in such a way
that $\varphi (T_j\cap U')\subset \{ \, (v_1,v_2)\in V\, |\, v_1<0\, \} $.
Let $\alpha(r),r\geqslant 0$, be real smooth monotone function such that
$\alpha(r) = 1$ for $r\in [0,\frac{4}{3}]$ and $\alpha(r) = 0$ for
$ r\geq \frac{5}{3}$. For $\overline Q =\sigma _{j_1}^{\overline \delta _1}
\cdot \dots \cdot \overline \sigma _{j_k}^{\delta _k}$
consider a smooth isotopy
$F_{\overline Q,t}:V\times E_R\to V\times E_R$, where
$$F_{\overline Q,t}=  F_{2k,\overline \sigma _{j_1}^{-\delta _1},t}\circ
\dots \circ
F_{k+2,\overline \sigma _{j_{k-1}}^{-\delta _{k-1}},t}\circ
F_{k+1,\overline \sigma _
{j_k}^{-\delta _k},t}\circ  F_{k,\overline \sigma _{j_k}^{\delta _k},t}
\circ \dots \circ F_{1,\sigma _{j_1}^{\delta _1},t}
$$
and $F_{l,\sigma _{j}^{\delta _j},t}$ is
given by functions
$$F_{l,\overline \sigma _{j}^{\delta _j},t}( s_1,s_2,y)=(s_1,s_2,
f_{l,\overline \sigma _{j}^{\delta _j},t}(s_1,s_2,y),$$
where
$$
\align
 & f_{l,\overline \sigma _{j}^{\delta _j},t}(s_1,s_2,y)= \\
= & \cases
 y\, ,  &  s_1\leq l-1 \, ; \\
 2j-2 +(y-2j+2)e^{i\pi \delta _j \alpha (s_2)\beta (s_1-l+1)\alpha(|y-2j+2|)t }
\, , & l-1\leq s_1\leq l\, ; \\
 2j-2 +(y-2j+2)e^{i\pi \delta _j \alpha (s_2)\alpha(|y-2j+2|)t}\, ,
& s_1\geq l\,
\endcases
\endalign
$$
if $l\leq k$, and
$$
\align
 & f_{l,\overline \sigma _{j}^{\delta _j},t}(s_1,s_2,y) = \\
= &
\cases
 y\, , &  s_1\geq l \, ; \\
 2j-2 +(y-2j+2)e^{i\pi \delta _j \alpha (s_2)\beta (l-s_1)\alpha(|y-2j+2|)t}
\, , & l-1\leq s_1\leq l\, ; \\
 2j-2 +(y-2j+2)e^{i\pi \delta _j \alpha (s_2)\alpha(|y-2j+2|)t}\, ,
& s_1\leq l-1\,
\endcases
\endalign
$$
if $l\geq k+1$.
One can check that
\roster
\item $F_{\overline Q,t}$ is the identity map over a neighbourhood of the
boundary of $V$ for all $t$;
\item $F_{\overline Q,1}(\overline{\Cal B})=\overline Q$, where
$\overline{\Cal B}=
\{ b_1(x(s_1,0))\equiv 1, \dots ,b_p(x(s_1,0)\equiv 2p-3 \} $,
$0\leq s_1\leq k$, is the trivial geometric braid.
\item $F_{\overline Q,1}(\Cal B)=\{ (s_1,s_2,-1),(s_1,s_2,1),\dots ,
(s_1,s_2,2p-3) \} $ are constant sections over $$\{ k-\frac{1}{3} < s_1<
k+\frac{1}{3}\} \times \{ 0<s_2<2 \} ;$$
\endroster

Such isotopy $F_{\overline Q,t}$ will be called {\it a
$\overline Q$-twisting-untwisting of constant sections
with support $\varphi ^{-1}(V)$ and with
center ($V_0, z_0$)}, where $V_0=\varphi ^{-1}(\{ k-\frac{1}{3} < s_1<
k+\frac{1}{3}\} \times \{ 0<s_2<2 \} $ and $z_0=\varphi ^{-1}((k,0))$.

Let $\widetilde F_{\overline Q,t}=\varphi ^*(F_{\overline Q,t})$. Denote
again by $B_1$ its image $\widetilde F_{\overline Q,1}(B_1)$.

In the notation of the definition of $g$-base $\Gamma $ and its tubular
neighbourhood, we change the $g$-base $\Gamma $ to an equivalent one
taking $z_0$ instead of $o,$ changing each path $T_j$ to a path starting
at $z_0$ and coinciding with $T_j$ outside the disc $E_{o}$.
We change $E_{o}$ to a disc $E_{z_0}\subset V_0$ with center at $z_0$ and
choose new neighbourhoods $U_{j}$ contained in the old neighbourhoods $U_j$.
In the sequel, we denote again by $o$ the point $z_0$.

By construction of the $\overline Q$-twisting-untwisting $F_{\overline Q,t}$,
the braid monodromy factorization of the curve $B_1$ will be
$$\Delta ^2 = \prod\limits_{j=1}^{n}Q^{-1}Q_j^{-1}\sigma _1^{\rho _j}Q_jQ,
$$
that is, the braid monodromy factorizations of the curves $B_1$ and $B_2$
coincide.

{\it Step V.} Let $u'_j$ be the point chosen in Step III, $u'_j\in T_j$,
and let $U'_{\Gamma }\subset U_{\Gamma }$ be a tubular neighbourhood of the
$g$-base $\Gamma $ such that $U'_{\Gamma }\cap U(\partial U_{\Gamma })=
\emptyset $, where $U(\partial U_{\Gamma })$ is a neighbourhood of the
boundary $\partial U_{\Gamma }$ of $U_{\Gamma }$.

Let us show that there exists a smooth isotopy $F_t:\F \to \F$ having
properties $(1)$ - $(3)$ of Theorem 7.1 and such that
\roster
\item $F_t$ is the identity map over the complement of $U_{\Gamma }$\, ;
\item $F_1(B_1)=B_2$ over $U'_{\Gamma }\subset U_{\Gamma }$.
\endroster

To show it, for each $j$ consider the geometric braids
$$\overline{ \Cal B}_1=(b_{1,j,1}(x),\dots ,b_{1,j,p}(x))\, \, \, \text{and}\,
\, \, \overline{ \Cal B}_2=(b_{2,j,1}(x),\dots ,b_{2,j,p}(x)),$$
where $x$ is moving along $T'_j\subset T_j$ starting at $o$ and ending at
$u'_j$. These geometric braids are representatives of elements $\beta _{i,j}
\in B[\Bbb C_o,K]$, $i=1,2$. Since the corresponding factors of the braid
monodromy factorizations for $B_1$ and $B_2$ coincide, then
$$\beta _{1,j}^{-1}H_1^{\nu _j}\beta _{1,j}=\beta _{2,j}^{-1}
H_1^{\nu _j}\beta _{2,j}=Q_{j}^{-1}H_1^{\nu _j}Q_{j}.$$
Thus, $\beta _{j}=\beta _{1,j}\beta _{2,j}^{-1}\in C(H_1^{\nu _j})$.
By Theorem 5.1 the element $\beta _{j}$ can be written in the form
$$\beta _{j}=\mu _{j,1}\dots \mu _{j,k_j},$$
where each $\mu _{j,i}$ coincides with either $H_r^{\delta _{j,i}}$, where
$\delta _{j,i}=\pm 1$ and $r=1,3,\dots , p$, or a full-twist $\Dl _{1,r}^
{2\delta _{j,i}}$, $r=3,\dots , p$, defined by the system of paths
$(\sigma _1,\dots ,\sigma _{p-1})$.

To each such $\beta _j$ we associate a ``twisting-untwisting" of $V_j\times
\Bbb C^1$ similar to the one described in Step IV. Namely, we consider again
$V_j=(-1,2k_j+1)\times (0,2)$, and for each $\mu _{j,i}$ we define a smooth
isotopy
$$F_{j,\mu _{j,i},t}:V_j\times \Bbb C\to V_j\times \Bbb C,$$
and associate to $\beta _{j}=\mu _{j,1}\dots \mu _{j,k_j}$ the composition
$$F_{j,\beta _{j},t}=F_{j,2k_j,\mu _{j,1}^{-1},t}\circ \dots \circ
F_{j,k_j+2,\mu _{j,k_j-1}^{-1},t}\circ F_{j,k_j+1,\mu _{j,k_j}^{-1},t}\circ
F_{j,k_j,\mu _{j,k_j},t}\circ \dots \circ
F_{j,1,\mu _{j,1},t}$$
as follows. If $\mu _{j,l}=H(\sigma _r)^{\delta _{j,l}}$, then
$F_{j,l,\mu _{j,l},t}=F_{l,\overline \sigma _r^{\delta _{j,l}},t}$ which
was defined in Step IV (the number $k$ in the definition of $F_{l,\overline
\sigma _{r}^{\delta _{j,l}},t}$ is equal to $k_j$ in our case). If
$\mu _{j,l}=\Dl _{1,r}^{2\delta _{j,l}}$, then $F_{j,l,\mu _{j,l},t}$ is
defined similarly, namely, it is given by
$$F_{j,l,\mu ^{\delta _{j,l}},t}( s_1,s_2,y)=(s_1,s_2,
f_{j,l,\mu ^{\delta _{j,l}},t}(s_1,s_2,y),$$
where
$$
\align
& f_{j,l,\mu ^{\delta _{j,l}},t}(s_1,s_2,y)= \\
= & \cases
 y\, ,  &  s_1\leq l-1 \, ; \\
 r-2 +(y-r+2)e^{2i\pi \delta _{j,l} \alpha (s_2)\beta (s_1-l+1)
 \gamma (|y-r+2|)t }
\, , & l-1\leq s_1\leq l\, ; \\
 r-2 +(y-r+2)e^{2i\pi \delta _{j,l} \alpha (s_2)\alpha(|y-r+2|)t}\, ,
& s_1\geq l\, ,
\endcases
\endalign
$$
if $l\leq k_j$, and
$$
\align
 & f_{l,\mu _{j,l}^{\delta _{j,l}},t}(s_1,s_2,y)= \\
= & \cases
 y\, , &  s_1\geq l \, ; \\
 r-2 +(y-r+2)e^{2i\pi \delta _{j,l} \alpha (s_2)\beta (l-s_1)
 \gamma(|y-r+2|)t}
\, , & l-1\leq s_1\leq l\, ; \\
 r-2 +(y-r+2)e^{2i\pi \delta _{j,l} \alpha (s_2)\gamma(|y-r+2|)t}\, ,
& s_1\leq l-1\,  ,
\endcases
\endalign
$$
if $l\geq k_j+1$, where $\alpha(s)$, $\beta(s)$, and $\gamma (s)$,
$s\geqslant 0$, are real smooth monotone functions such that $\alpha(s) = 1$
for $s\in [0,\frac{4}{3}]$ and $\alpha(s) = 0$ for $ s\geq \frac{5}{3}$,
$\beta (s)=0$ for $s\in [0,\frac{1}{3}]$ and $\beta (s)=1$ for $s\geq
\frac{2}{3}$, and $\gamma(s) = 1$ for $s\in [0,r-1]$ and $\gamma(s) = 0$ for
$ s\geq r-\frac{1}{2}$.

Let us choose a neighbourhood $W_j\subset E_j$ containing the part of diameter
$d'_j$ connecting the points $u'_j$ and $u^{\prime \prime}_j$, and such that
there exists a diffeomorphism $\phi _j: W_j\to V_j$ such that $\phi _j(d'_j)=
\{\, (s_1,s_2)\in V_j\, \, |\, \, 0\leq s_1\leq 2k_j, \, s_2=0\, \}$,
$\phi _j(u_j)=(0,0)$. The diffeomorphism $\phi _j$ and the isotopy
$F_{j,\beta _{j},t}$ allow us to define a smooth isotopy $\widetilde
F_{j,\beta _{j},t}=(\phi _j^{-1}\times Id)\circ F_{j,\beta _{j},t}\circ
(\phi _j\times Id):W_j\times \Bbb C^1\to W_j\times \Bbb C^1$ which can be
extended to a smooth isotopy such that $\widetilde F_{j,\beta _{j},t}$ is
the identity map outside $W_j\times \Bbb C^1$. Let $F_t$ be the
composition of the constructed isotopies $\widetilde F_{j,\beta _{j},t}$,
$j=1,\dots ,n$. Denote again by $B_1$ its image $F_1(B_1)$ and by $E_j$ a
disc contained in $V_{0,j}$, where $V_{0,j}$ is the center of the
``twisting-untwisting"
$\widetilde F_{j,\beta _{j},t}$, which is defined verbatim   in Step IV.
We choose a new point in $E_j\cap T_j$ and denote it again by $u'_j$. By
construction of the isotopies $\widetilde F_{j,\beta _{j},t}$, for each $j$
the geometric braids
$$\overline{ \Cal B}_1=(b_{1,j,1}(x),\dots ,b_{1,j,p}(x))\, \, \, \text{and}\,
\, \,
\overline{ \Cal B}_2=(b_{2,j,1}(x),\dots ,b_{2,j,p}(x)),$$
where $x$ is moved along $T'_j\subset T_j$ starting at $o$ and ending at
$u'_j$, are two representatives of the same element of $B_p[\Bbb C^1_o,K]$.
Therefore by Lemma 6.10, there exists a smooth isotopy $F_t:\F \to \F$
having properties $(1)$ - $(3)$ of Theorem 7.1 such that
\roster
\item $F_t$ is the identity map over the complement of the union of
small neighbourhoods $U_j$ of the paths $T'_j$\, ;
\item $F_t$ is the identity map over $(\cup E_j)\cup E_{o}$\, ;
\item $F_1(B_1)=B_2$ over some tubular neighbourhood
$U'_{\Gamma }\subset U_{\Gamma }$ of the $g$-base $\Gamma $\,.
\endroster

Denote again by $B_1$ its image $F_{1}(B_1)$.
The obtained curves $B_1$ and $B_2$ coincide over the tubular neighbourhood
$U'_{\Gamma }$ of the $g$-base $\Gamma $.

{\it Step VI.} The complement
$\Bbb P ^1 \setminus U'_{\Gamma }$ is simply connected. Let $U_{\infty }$ be
a simply connected neighbourhood of $\Bbb P ^1 \setminus U'_{\Gamma }$ such
that $U_{\infty }$ is diffeomorphic to a disc and such that $u_j\not\in
U_{\infty }$ for all $j=1,\dots ,\# N$. Denote by $V=\pi ^{-1}(U_{\infty})$.
Then $\pi : V\to U_{\infty }$ is a trivial fibering with fibres $\pi ^{-1}(x)
\simeq \Bbb P ^1$ and $L\cap V$ is a section. Put $V_0= V\setminus L$ and
${\Cal B}_i=V_0\cap B_i$, $i=1,2$. Then $V_0\simeq U_{\infty }\times \C ^1$
and $\pi _{|V_0}$ coincides with projection on the second factor. We apply
Remark 6.7 and Lemma 6.9 to ${\Cal B}_1$ and ${\Cal B}_2$ to obtain
a smooth isotopy $F_t$ having all properties (1) - (4) of
Theorem 7.1.
\enddemo
\bk

\subheading{\S 8.\ Equivalence of braid monodromy factorizations  and
diffeomorphism types of surfaces}

In this Section we prove Theorem 2.

Let $F :(\CP ^2, B_1)\to (\CP ^2, B_2)$ be a diffeomorphism of pairs
having the properties described in Corollary 7.2. This diffeomorphism induces
an isomorphism $F^{*} :\pi _1(\CP ^2\setminus B_2)\to \pi _1
(\CP ^2\setminus B_1)$.

By Proposition 1 in \cite{Kul1}, the set of non-equivalent generic morphisms
of degree $N$ with discriminant curve $B\subset \CP ^2$ is in one-to-one
correspondence with the set of epimorphisms from $\pi _1 (\CP ^2\setminus B)$
to the symmetric group $\Sigma _N$ satisfying some additional conditions
(see details in \cite{Kul1}). Since Chisini's Conjecture holds for
$B_1\subset \CP ^2$, then there exists such a unique   epimorphism from
$\pi _1 (\CP ^2\setminus B_1)$, which must coincide with the epimorphism
$f_{1*}:\pi _1 (\CP ^2\setminus B_1)\to \Sigma _N$ induced by $f_1$, where
$N=\deg f_1$. Therefore, for $B_2$, there exists such a unique
epimorphism, which must
coincide with $f_{1*}\circ F^{*}=f_{2*}:\pi _1 (\CP ^2\setminus B_2)\to
\Sigma _N$.
Consequently, the diffeomorphism $F :\CP ^2\setminus B_1\to \CP ^2\setminus
B_2$ can be
lifted to a diffeomorphism $\Psi _0:S_1
\setminus f_1^{-1}(B_1)\to S_2\setminus f_2^{-1}(B_2)$.

In \cite{Kul2}, one can find a method how to reconstruct a surface $S$ and a
finite morphism $f:S\to \CP ^2$ branched along $B\subset \CP ^2$ if we know
the homomorphism $f_{*}:\pi _1(\CP ^2\setminus B)\to \Sigma _N$. This method
is based on the presentation $S$ as $N$ copies of $\CP ^2$ with ``standard
cuts" pasted together along these cuts (to do such pasting together, we use
the geometric desciption of the finite presentation $\pi _1(\CP ^2 \setminus
B)$ in terms of ``shadows" and ``screens" described in \cite{Kul3}).
Using this method, it is easy to see that the diffeomorphism $\Psi_0$ is
uniquely extended to a homeomorphism $\Psi:S_1\to S_2$.

Let $U\subset S_1$ be a neighbourhood of an ordinary cusp of $B_1$ such that
$F_{|U}$ is holomorphic. It is well-known that if $f:X\to U$ is a
three-sheeted covering of $U=\{ \, (x,y)\in \C ^2\, \mid \, |x|<1,\, |y|<1\,
\} $ branched along a curve given by $y^2=x^3$ and such that $f$ is not
Galois covering, then such $f$ is unique. Therefore, the homeomorphism
$\Psi $ is holomorphic (in particular, $\Psi $ is smooth) in $f_1^{-1}(U)$.
Similarly, $\Psi $ is smooth in $f_1^{-1}(U)$, where $U$ is a neighbourhood
of a node of $B_1$ or a tangent point of $B_1$ and a fiber of the
projection $\pi $.

Let $z\in B_1$ be a non-singular point and let
$$U_1\simeq \{ \,
(x_1,y_1)\in \C ^2\, \mid \, |x_1|<1,\, |y_1|<1\, \} $$
be a neighbourhood
of $z$ in $\CP ^2$, where $(x_1,y_1)$ local holomorphic coordinates in
$\CP ^2$ such that $y_1=0$ is a local equation of $B_1$ and the projection
$\pi $ is given in $U_1$ by $(x_1,y_1)\mapsto x_1$. Similarly, let $U_2=
F(U_1)$ be a neighbourhood of $F(z)$ and let $(x_2,y_2)$ be local holomorphic
coordinates in $U_2$ such that $y_2=0$ is a local equation of $B_2$ and the
projection $\pi $ is given in $U_2$ by $(x_2,y_2)\mapsto x_2$. We have
$x_2=g_1(x_1)$ and $y_2=g_2(x_1,y_1)$, where $g_1$ and $g_2$ are smooth
functions and $g_2$ is holomorphic in $y_1$. Therefore $g_2$ can be writen in
the form
$$g_2(x_1,y_1)=\sum\limits_{n=1}^{\infty }a_n(x_1)y_1^n, $$
where all $a_n(x_1)$ are smooth and $a_1(x_1)\neq 0$ in $U_1$.

Each preimage $f_1^{-1}(U_1)$ and $f_2^{-1}(U_2)$ consists of
$N-1$ connected components $U_{1,1}$,..., $U_{1,N-1}$ and
$U_{2,1},\dots ,\, U_{2,N-1}$, respectively. Let $f_1$ (resp. $f_2$) is
non-ramified in $\cup_{j=2}^{N-1}U_{1,j}$ (resp. in $\cup_{j=2}^{N-1}U_{2,j}$).
Therefore, $\Psi $ is smooth in $\cup_{j=2}^{N-1}U_{1,j}$.
Besides, there exist local holomorphic cordinates $(u_1,v_1)$ in $U_{1,1}$
(resp. $(u_2,v_2)$ in $U_{2,1}$) such that $f_1$ is given in $U_{1,1}$
(resp. $f_2$ in $U_{2,1}$) by $y_1=u_1^{2}$, $x_1=v_1$ (resp.
$y_2=u_2^{2}$, $x_2=v_2$). Consequently, $\Psi $ is given by
$$\align
u_2 & = u_1 (\sum\limits_{n=1}^{\infty }a_n(v_1)u_1^{2n-2})^{\frac{1}{2}}
\, ; \\
v_2 & = v_1 \, .
\endalign
$$
It is easy to see that
$$(\sum\limits_{n=1}^{\infty }a_n(v_1)u_1^{2n-2})^{\frac{1}{2}}  $$
is a smooth function, since all $a_n(v_1)$ are smooth and $a_1(v_1)\neq 0$.
\bk

\end

\end